\input amstex
\input epsf
\documentstyle{amsppt}
\magnification 1200
\nologo
\def\epsfsize#1#2{\hsize}
\TagsOnRight
\voffset -2pc


\def\on#1{\expandafter\def\csname#1\endcsname{{\operatorname{#1}}}}
\on{lf} \on{id} \on{dist} \on{deg} \on{Hom} \on{sign} \let\Tl\widetilde
\def\fn#1{\expandafter\def\csname#1\endcsname{\operatorname{#1}\def\shortcut
{\if(\next{}\else\,\fi}\futurelet\next\shortcut}} \fn{coker} \fn{rel}
\let\inv\leftrightarrow \let\emb\hookrightarrow \let\imm\looparrowright
\def\R{\Bbb R} \def\?{@!@!@!@!@!} \def\RP{\R\? P} \def\z#1#2{#1\?#2\?}
\def\Z{{\Bbb Z}\def\shortcut{\if/\next{}\z\fi}\futurelet\next\shortcut}
\let\tl\tilde \redefine\o{\circ} \let\Cap\smallfrown \let\but\setminus
\let\x\times \let\eps\varepsilon \let\phi\varphi \def\T{_{\?\text{\fiverm T}}}
\def\N{\Bbb N}  \def\Cl#1{\overline{#1}}

\topmatter 
\thanks Partially supported by the Russian Foundation for Basic Research
Grant No\. 02-01-00014. \endthanks
\address Steklov Mathematical Institute, Division of Geometry and Topology;
\newline $\Gamma$C$\Pi$-1, ul. Gubkina 8, Moscow 117966, Russia \endaddress
\curraddr University of Florida, Department of Mathematics;\newline
358 Little Hall, PO Box 118105, Gainesville, FL 32611-8105, U.S. \endcurraddr
\email melikhov\@mi.ras.ru, melikhov\@math.ufl.edu
\endemail

\title Sphere eversions and realization of mappings\endtitle
\author Sergey A. Melikhov \endauthor

\abstract\nofrills
P. M. Akhmetiev used a controlled version of the stable Hopf invariant
to show that any (continuous) map $N\to M$ between stably parallelizable
compact $n$-manifolds, $n\ne 1,2,3,7$, is realizable in $\R^{2n}$, i.e\.
the composition of $f$ with an embedding $M\i\R^{2n}$ is $C^0$-approximable
by embeddings.
It has been long believed that any degree $2$ map $S^3\to S^3$, obtained by
capping off at infinity a time-symmetric (e.g\. Shapiro's) sphere eversion
$S^2\x I\to\R^3$, was non-realizable in $\R^6$.
We show that there exists a self-map of the Poincar\'e homology $3$-sphere,
non-realizable in $\R^6$, but every self-map of $S^n$ is realizable in
$\R^{2n}$ for each $n>2$.
The latter together with a ten-line proof for $n=2$, due essentially to
M. Yamamoto, implies that every inverse limit of $n$-spheres embeds in
$\R^{2n}$ for $n>1$, which settles R. ~Daverman's 1990 problem.
If $M$ is a closed orientable $3$-manifold, we show that there exists
a map $S^3\to M$, non-realizable in $\R^6$, if and only if $\pi_1(M)$
is finite and has even order.

As a byproduct, an element of the stable stem $\Pi_3$ with non-trivial
stable Hopf invariant is represented by a particularly simple immersion
$S^3\imm\R^4$, namely the composition of the universal $8$-covering over
$Q^3=S^3/\{\pm1,\pm i,\pm j,\pm k\}$ and an explicit embedding
$Q^3\emb\R^4$.
\endabstract
\endtopmatter

\document
Let $N^n$ and $M^{2n-k}$, $0\le k\le n$, be orientable smooth manifolds,
where $N$ is compact, and let $f\:N\to M$ be a self-transverse smooth map
(i.e\. for every $x\ne y$ with $f(x)=f(y)$ the tangent space $T_{f(x)}M$
is generated by $df(T_xN)$ and $df(T_yN)$).
Consider the $k$-dimensional manifold $\Sigma_f/\tau$ of double points
of $f$, that is, unordered pairs $\{x,y\}$ of distinct points of $N$ such
that $f(x)=f(y)$, and let $\eta$ denote the line bundle, associated with
the $2$-cover $\Sigma_f\to\Sigma_f/\tau$ by the ordered pairs $(x,y)$.
Thus $\Sigma_f$ is the orientable submanifold
$\{(x,y)\mid x\ne y, f(x)=f(y)\}$ of the deleted product
$\tilde N:=N\x N\but\Delta_N$, which can be thought of as the preimage of
$\Delta_M$ under $f^2\:\tl N\to M\x M$; finally, $\tau$ is the factor
exchanging involution $(x,y)\inv (y,x)$ on $\tl N$.
We define the {\it stable Hopf invariant} $h(f)$ to be
$(w_1(\eta))^k[\Sigma_f/\tau]\in\Z/2$.
Note that in the case $k=0$ this is just the parity of the number of double
points of $f$.
In general, $h(f)$ can be alternatively described as the $\bmod 2$ degree
(summed over the components) of a map $\Sigma_f/\tau\to\RP^k$,
classifying the line bundle $\eta$.

If $M=\R^{2n-k}$ and $f$ happens to be an immersion with trivial normal
bundle, $h(f)$ coincides with the stable Hopf invariant $h([f,\Xi])$ of
the element of the $n$-th stable stem $\Pi_n=\pi_{n+m}(S^m)$, $m\ge n+2$,
represented by $f$ with an arbitrary framing $\Xi$ (Akhmetiev--Sz\H ucs
\cite{AS}; see also Koschorke--Sanderson \cite{KS; p\. 203} for $k=0$ and
Lemma 6(a) below).
We recall that $h([f,\Xi])$ is, by definition, the $\bmod 2$ reduction of
the (usual integral) Hopf invariant of any element of $\pi_{2n+1}(S^{n+1})$
whose suspension is $[f,\Xi]$.
In other words, $h([f,\Xi])$ is the parity of the linking number between
$\bar f(N)$ and $\bar f^+(N)$, where $\bar f\:N^n\to\R^{2n+1}$ is any
embedding projecting to $f$ along the last $k+1$ coordinates, and
$\bar f^+$ its pushoff along the last vector of the framing $\bar\Xi$,
obtained by adjoining $k+1$ constant vectors to $\Xi$.

By a classical theorem of J. F. Adams (see \cite{A4}, where a novel
immersion-theoretic approach is developed), $h(f)=0$ if $f$ is an immersion
into $\R^{2n-k}$ with trivial normal bundle and $n\ne 1,3,7$.
The figure $8$ immersion $S^1\imm\R^2$ with one double point clearly has
$h(f)\ne 0$.
Max and Banchoff \cite{MB} showed that any eversion of $S^2$, regarded as
a level-preserving immersion $\phi\:S^2\x I\imm\R^3\x I$, has an odd number
of quadruple points in general position (see also \cite{Hu}, \cite{N}).
By a theorem of M. Freedman \cite{Fre; Lemma 2} (see also
\cite{K; Theorem F(b)}) this means that after capping off at infinity to
an immersion $\hat\phi\:S^3\imm\R^4$, we obtain a representative for
an element of the stable stem $\Pi_3$ with nontrivial $h([\hat\phi])$.
(Other such explicit representatives are given in \cite{K}, \cite{Ca},
\cite{N} and in Example 4 below.)
The Max--Banchoff and Freedman theorems do not extend for $n=7$, for no
orientable $7$-manifold generically immerses in $\R^8$ with an odd number
of $8$-tuple points \cite{Ec}.

It is implicit in the work of P. M. Akhmetiev \cite{A1}, \cite{A2} that,
if $n>2$, the {\it controlled stable Hopf invariant}
$H(f):=(w_1(\eta))^k\in H_0(\Sigma_f/\tau;\Z/2)$ of the self-transverse map
$f\:N^n\to M^{2n-k}$ vanishes if and only if the composition of $f$ with
$M\subset M\times\R^k$ is $C^0$-approximable by embeddings (see Theorem 5
below).
The immersion-theoretic proof of the easy part of the Adams theorem \cite{AS}
(see also \cite{Mi; Corollary 1}) shows that $H(f)=0$ for any self-transverse
map $f\:N^n\to M^{2n-k}$ between stably parallelizable manifolds, where $k$
is not less than the maximal power of $2$ dividing $n+1$.
On the other hand, by examples of Miller and Wells (see \cite{Mi; Theorem 3})
for every $k$ there exist infinitely many $n$ such that $h(f)\ne 0$ (hence
also $H(f)\ne 0$) for some self-transverse immersion $f\:S^n\imm\R^{2n-k}$
(with possibly nontrivial normal bundle).

\example{\vskip-8pt Example 1} Any self-transverse immersion
$S^n\overset f\to\imm\R^{n+1}\overset i\to\i\R^{2n}$ with $h(f)\ne 0$
(such as a capped off eversion of $S^2$ in $\R^3$) is not approximable by
embeddings.
Indeed, any self-transverse immersion $\bar f\:S^n\imm\R^{2n}$, projecting
to $f$ along the last $n-1$ coordinates, has the same stable Hopf invariant,
hence an odd number of double points.
Pick such an $\bar f$, $\eps$-close to $if$, and suppose that there exists
an embedding $g\:S^n\emb\R^{2n}$, $\eps$-close to $if$ for a sufficiently
small $\eps>0$ (defined below).
Since $\Sigma_{\bar f}\i\Sigma_f$, some connected component $C$ of
$\Sigma_f/\tau$ contains an odd number of the points of
$\Sigma_{\bar f}/\tau$.
On the other hand, for any map $h\:S^n\to\R^{2n}$, $\eps$-close to $f$,
the set $\Sigma_h$ lies in the $\eps$-neighborhood $O_\eps$ of
$\Sigma_f\cup\Delta_{S^n}$.
Hence any self-transverse $\eps$-homotopy $H\:S^n\x I\to\R^{2n}\x I$ between
$\bar f$ and $g$ defines a $\tau$-equivariant null-bordism
$\Sigma_H\i O_\eps\x\Delta_I\i\widetilde{S^n\x I}$ for $\Sigma_{\bar f}$.
If $\eps$ is small enough, the $\eps$-neighborhood of $C$ is disjoint from
the $\eps$-neighborhoods of the other components and of $\Delta_{S^n}$, so
the set $C\cap\Sigma_{\bar f}/\tau$ of odd cardinality is null-bordant. \qed
\endexample

The manifold $M^{2n-k}\x\R^k$ appears to be the simplest reasonable target
space for controlled embedding theory, because every map
$N^n\to M^{2n-k}\x\R^{k+1}$ is approximable by embeddings by general
position.
On the other hand, the critical value of the parameter $k$ is $n$, because
every map $f\:N^n\to M^{n-1}\i M^{n-1}\x\R^{n+1}$ is approximable by
embeddings, provided that $N$ is unsophisticated enough to embed in $\R^{n+1}$
(say, if $N$ any homology $n$-sphere%
\footnote{Every homology sphere bounds a contractible topological manifold
\cite{Ke; Corollary}, \cite{FQ} (for $n\ge 5$ see also \cite{AG}), whose
double has to be $S^{n+1}$ by Seifert--van Kampen and the generalized
Poncar\'e Conjecture.
This yields an approximation of $f$ by topological embeddings, which in turn
can be approximated by smooth ones \cite{Ha}.}%
).
Looking forward to getting rid of the target manifold $M^n\x\R^n$ by
embedding it into the Euclidean space of the same dimension, we arrive at
the requirement that $M$ be stably parallelizable.

\definition{Definition}
We say that a (continuous) map $f\:N\to M$ between smooth $n$-manifolds,
where $N$ is compact and $M$ is stably parallelizable, is {\it realizable}
(by embeddings) in $\R^{2n}$, if the composition of $f$ and some smooth
embedding $i\:M\emb\R^{2n}$ is $C^0$ $\eps$-close to a smooth embedding
for each $\eps>0$.
The choice of $i$ does not matter, since its normal bundle is always trivial
(this is clearly so outside a point \cite{KM; Lemma 3.5}, thus it remains to
take into account that an arbitrary embedding $S^n\emb\R^{2n}$ has trivial
normal bundle \cite{KM; Lemma 8.3}).
For $n>2$ we could even do with topological embeddings, since they are
approximable by smooth ones \cite{Ha}.
Note that if $g\:L\to M$ and $f\:M\to N$ are realizable in $\R^{2n}$, so is
$fg$.
\enddefinition

\example{Example 2} Every map $f\:S^n\to\R^n$ is realizable in $\R^{2n}$.
To show this, we may assume without loss of generality that $f$ is
self-transverse.%
\footnote{A smooth map with compact domain is $C^\infty$-approximable by
self-transverse ones \cite{GG; III.3.2}.}
Then a generic point in the topological frontier of $f(S^n)$ has exactly
one preimage.
We will show that every map $f\:S^n\to S^n$ such that some $p\in S^n$ has
exactly one preimage $q$, is realizable in $\R^{2n}$.
Indeed, let $h\:S^n\but\{q\}\to\R^n$ be a homeomorphism; then
$g\:S^n\to S^n\x\R^n$, defined by $g(x)=(f(x),\dist(f(x),p)h(x))$ for
$x\ne q$ and $g(q)=(p,0)$, is an embedding. \qed
\endexample

P. M. Akhmetiev proved \cite{A1}, \cite{A2} (see also \S3 and \cite{A3} for
alternative arguments) that if $n\ne 1,2,3,7$, every map $f\:N\to M$ between
stably parallelizable $n$-manifolds, where $N$ is compact, is realizable in
$\R^{2n}$.
Indeed, it follows from the Adams theorem that $H(f)=0$ if $n\ne 1,3,7$
(see Lemma 9 below).
It is not hard to see that a map $f\:S^1\to S^1$ is realizable in $\R^2$ iff
$\deg(f)\ne 0,1,-1$ \cite{Mc}, \cite{Si}.
It was shown in \cite{Sk; \S5} that the composition of the $2$-cover
$S^n\x I\to\RP^n\x I$ and an immersion $\RP^n\x I\imm\R^{n+1}\i\R^{2n+1}$
(existing for $n=1,3,7$) is not approximable by embeddings.

An interesting attempt to find maps $S^n\to S^n$, $n=3,7$, non-realizable in
$\R^{2n}$, occurred in \cite{A1}, \cite{A2}.
Any sphere eversion $\phi\:S^{n-1}\x I\to\R^n$ (which exists iff $n=1,3$ or
$7$) can be ``capped off'' at infinity so as to produce a degree two map
$f_\phi\:S^n\to S^n$.
A {\it time-symmetric} eversion passes transversely at $t=1/2$ through the
double cover of some immersion $\RP^{n-1}\imm\R^n$ and returns back by the
same route, precomposed with the antipodal involution of $S^{n-1}$;
time-symmetric eversions also exist for $n=1,3$ and $7$, see \cite{A2}.
It was claimed in \cite{A1}, \cite{A2} that if $\phi$ is time-symmetric,
the map $f_\phi$ (which in this case factors through the $2$-cover
$S^n\to\RP^n$) is non-realizable in $\R^{2n}$.
That this is not so first became clear from Example 7 in \S1, where
the degree two map $S^3\to S^3$ yielded by Shapiro's original time-symmetric
eversion \cite{Fra} is shown to realize in $\R^6$.
With only slightly increasing difficulty, we further show in \S1 that any
capped off time-symmetric eversion of $S^{n-1}$ in $\R^n$, $n=3,7$, is
realizable in $\R^{2n}$.

Our starting point is the following

\proclaim{Theorem 1} Let $N$ be a $\Z/2$-homology $n$-sphere, $M$ a stably
parallelizable $n$-manifold, $n>2$, and $f\:N\to M$ a self-transverse map.
Then $f$ is realizable in $\R^{2n}$ iff the restriction of the projection
$p_1\:N\x N\to N$ to each $\tau$-invariant connected component of
$\Sigma_f$ is of even degree.%
\footnote{We follow the convention that
the degree of a map of a non-compact $n$-manifold to a compact $n$-manifold
is zero (rather than undefined).}
\endproclaim

Note that the analogous assertion for $n=1$ is not true, for there are no
$\tau$-invariant components in $\Sigma_f$, where $f\:S^1\to S^1$ is the
$3$-cover, or more generally any odd degree map.

The proof, which is more or less in the spirit of Akhmetiev's arguments in
\cite{A1}, \cite{A2}, \cite{A3} (and hopefully clarifies them somewhat), is
postponed until \S3.
Note that the condition in Theorem 1 trivially holds if every component of
$\Sigma_f$ is either non-compact or mapped by $p_1$ onto a proper subset of
$N$.
The reader may find it exciting to verify by inspection that this is
the case for the degree two map $f\:S^3\to S^3$ obtained by capping off
Morin's eversion \cite{Fra}.
We will eventually show in \S1 that the same is true of any eversion.

\example{Example 3} The $2$-cover $f\:S^3\to\RP^3$ is not realizable in
$\R^6$, for $\Sigma_f$ coincides with the {\it antidiagonal}
$\nabla_{S^3}=\{(x,-x)\mid x\in S^3\}$ which projects homeomorphically onto
each factor of $S^3\x S^3$.
For the $4$-cover $f\:S^3\to S^3/(\Z/4)=L(4,1)$, the only $\tau$-invariant
component of $\Sigma_f$ is again $\nabla_{S^3}$, thus $f$ does not realize
in $\R^6$.
More generally, if $G$ is a finite group acting freely on $S^3$, the cover
$f\:S^3\to S^3/G$ is regular, hence the components of $\Sigma_f$ project
homeomorphically onto each factor of $S^3\x S^3$, and are in bijective
correspondence with the elements of $G\but\{1\}$, where elements of order
$2$ correspond to $\tau$-invariant components.
Thus, $f$ realizes in $\R^6$ iff $G$ has an odd order (using the Sylow
theorem and that orientable $3$-manifolds are parallelizable).
In particular, the $8$-cover of the quaternion space
$Q^3=S^3/\{\pm1,\pm i,\pm j,\pm k\}$ and
the $120$-cover of the Poincar\'e homology sphere $P^3=SO(3)/A_5$ do not
realize in $\R^6$.
Similarly, the $2$-cover $S^7\to\RP^7$ is not realizable in $\R^{14}$;
however, the $7$-dimensional lens space $S^7/(\Z/4)$ is not stably
parallelizable, see \cite{A4}.
\endexample

\remark{Remark} Using a generalization of Theorem 1 from \S3, we may
analogously conclude that for every odd $n$ the composition of the cover
$f\:S^n\to\RP^n$ and the inclusion $i\:\RP^n\i\RP^n\x\R^n$ is not approximable
by embeddings.
An elementary approach, valid in all dimensions, is as follows: let
$j\:S^n\i\R^{n+1}$ denote the inclusion and $p\:\R^{n+1}\to\R^n$
the projection; the map $f\x(pj)\:S^n\to\RP^n\x\R^n$ is an immersion with
one transverse double point; since $\Sigma_f=\nabla_{S^n}$, the argument
of Example 1 now applies.
In particular, this yields a proof of the Borsuk--Ulam theorem (stating that
there is no map $\phi\:S^n\to\R^n$ such that $f\x\phi\:S^n\to\RP^n\x\R^n$ is
an embedding).
Actually, it is well-known that $if$ is not even homotopic to an embedding
\cite{R}, \cite{DH}.
\endremark

\example{Example 4} By Example 3, for a finite group $G$ acting freely on
$S^3$, the cover $p\:S^3\to S^3/G$ has $h(p)\ne 0$ iff $G$ contains an odd
number of elements of order $2$.
In particular, $h(q)\ne 0$ for the universal $8$-cover $q\:S^3\to Q^3$ of
the quaternion space.
If $e\:Q^3\emb\R^4$ is a smooth embedding, it follows from Lemma 6 below that
$h(\phi)\ne 0$ for a self-transverse immersion $\phi$ in the regular homotopy
class of $eq$.
Thus $eq$ itself is a representative of an element of $\Pi_3$ with
non-trivial stable Hopf invariant.
There exists a particularly simple embedding $e\:Q^4\emb S^4$, apparently
referred to by E. Rees in \cite{R}, such that the closure of each
complementary domain is a tubular neighborhood of an embedded $\RP^2$
(compare \cite{Ma}).

Let $K\cup K'\i S^3$ be the $(2,4)$-torus link, i.e.\ the orbit of a pair
of points $(p,q,\frac12)$, $(p,-q,\frac12)$ under the action of $S^1$ on
$S^3=S^1*S^1$ given by $z\cdot(x,y,t)=(z+x,2z+y,t)$.
The orbit $\mu$ of the segment $(p,q)\x[\frac12,1)\cup\emptyset*\{q\}$ is
a M\"obius band, bounded by $K$ in the complement of $K'$, so that
its median circle $\emptyset*S^1$ is Hopf-linked with $K'$.
On the other hand, $K$ and $K'$ are interchanged by the
orientation-preserving involution $T\:(x,y,t)\mapsto(x,-y,t)$ of $S^3$.

Each of $K$ and $K'$ is unknotted in $S^3$, hence they bound disjoint
unknotted disks $D$, $D'$ in the upper and lower hemispheres of $S^4$,
separated by the equator $S^3$.
We may now subtract an embedded $2$-handle $H'$ from the lower hemisphere
$B^4$ along $D'$, so that the resulting manifold $M_0=\Cl{B^4\but H'}$ is
homeomorphic with $S^1\x D^3$, and then add an embedded $2$-handle $H$ to
$M_0$ along $D$.
The resulting $4$-manifold $M=M_0\cup H$ is easily seen to collapse onto
$\mu\cup D\cong\RP^2$, hence by uniqueness of regular neighborhood is
homeomorphic to a bundle over $\RP^2$ with fiber $D^2$.
On the other hand, its boundary $\partial M$ is the result of $0$-surgery
on the link $K\cup K'$, and the involution $T$ extends to an involution of
$S^4$ sending $M$ to its exterior $\Cl{S^4\but M}$.

It remains to construct a homeomorphism $\partial M\cong Q^3$.
The orbit of the segment $(p,-q)\x(0,\frac12]\cup\{p\}*\emptyset$ is an
annulus in the complement of $K$, cobounding $K'$ and $S^1*\emptyset$,
hence the link $K\cup K'$ is equivalent to $K\cup S^1*\emptyset$.
Consequently the pair $(\partial M_0,K)$ is homeomorphic to $(S^2\x S^1,K)$,
where $K$ is contained in the equatorial $S^1\x S^1\i S^2\x S^1$ in
the same way as it was in $S^1\x S^1\x\frac 12\i S^1*S^1$.
So $(S^2\x S^1,K)$ is homeomorphic to the mapping torus of the linear
homeomorphism of $(S^0*S^0*S^0,\,S^0*\emptyset*\emptyset)$, which is
antipodal on the first two factors and identical on the third.
The latter is homeomorphic to the quotient of $(D^2\x S^1,\{0\}\x S^1)$ by
an involution $\tau$ of the boundary that takes the meridian $y=0$ onto
the opposite meridian $y=\pi$ with the opposite orientation, and restricts
to the antipodal involution of the $(1,2)$-curve $y=2x$ (which can be
identified with the median circle of the M\"obius band $\mu$).
It is not hard to see that on the universal cover this involution is
the composition of the reflection through the line $y=2x$ along the $y$ axis
and the shift by $\pi$ along the $x$ axis.
Finally, $\partial M$ is obtained from $\partial M_0\cong D^2\x S^1/\tau$
by means of the $0$-surgery on $K=\{0\}\x S^1$, thus
$\partial M=D^2\x S^1/\tau'$, where $\tau'$ is $\tau$ conjugated by
the interchange $(x,y)\mapsto (y,x)$ of the meridian and the longitude.

On the other hand, the manifold $Q^3$ can be obtained by glueing together
the opposite faces of the cube by a clockwise rotation by $\pi/2$.
(The cube is a fundamental domain of the action of the quaternion group
of order $8$ on $S^3$, and its faces, edges and vertices are respectively
nonvoid intersections of the cube with one, two or three other fundamental
domains.)
After gluing together a pair of opposite faces (say, the top and the bottom),
$Q$ becomes the quotient of the solid torus by an involution of the boundary
that interchanges the meridian $y=0$ (represented by the horizontal edges in
the cube) with the $(4,1)$-curve $x=4y$ (made up of the vertical edges),
preserving their orientations in the direction of the $y$ axis.
On the universal cover the involution is easily seen to be the composition
of the reflection through the line $x=2y$ along the $x$ axis and the shift
by $\pi$ along the $y$ axis. \qed
\endexample

\remark{Remark} No such example can be constructed in dimension $7$.
Indeed, if $G$ is a finite group of even order acting freely on
a $\Z/2$-homology sphere $N$, by Milnor's and Smith's theorems
(see \cite{Bre}) there is precisely one element $g$ of order $2$ in $G$.

If $|G|$ is divisible by $4$, this and the Sylow theorem imply that $G$
contains an element $h$ of order $4$.
Since $N$ is a $\Z/2$-homology $7$-sphere, the Smith sequence shows that
$N/\left<g\right>$ has the same $\bmod 2$ cohomology ring as $\RP^7$.
Therefore, a map $N/\left<g\right>\to\RP^7$ classifying the $2$-cover
$N\to N/\left<g\right>$ is of an odd degree, and hence so is a map
$N/\left<h\right>\to S^7/(\Z/4)$ classifying the $4$-cover
$N\to N/\left<h\right>$.
On the other hand, if $N/G$ embeds in $S^8$, it is stably parallelizable.
However, there exist no odd degree maps of a stably parallelizable
$7$-manifold onto $S^7/(\Z/4)$ \cite{A4}.

If $|G|\equiv 2\pmod 4$, the Smith sequence for the $2$-cover
$N\to N/\left<g\right>$ shows that $N/G$ has the same $\Z/2$-cohomology
and same $\Z_{(2)}$-homology as $\RP^7$.
(Indeed, since $N$ is a $\Z/2$-homology sphere, its integral cohomology
is odd torsion in dimensions $1,\dots,6$, hence $N$ is a $\Z_{(2)}$-homology
sphere as well.)
Suppose that $S^8=A\cup B$ with $A\cap B\cong N/G$ (compare \cite{R}).
By Mayer--Vietoris, $H_i(N/G)\simeq H_i(A)\oplus H_i(B)$ for $1\le i\le 6$.
Say, $H_3(A;\Z_{(2)})\simeq\Z/2$ and $H_3(B;\Z_{(2)})=0$.
By the Bockstein sequence, $H_4(A;\Z/2)\ne 0$.
Hence $H^4(A;\Z/2)\ne 0$, contradicting the Alexander duality.
\endremark

\definition{Definition} For any map $f\:X\to Y$ the {\it double point map}
$f^{(2)}\:\Sigma_f/\tau\to Y$ of $f$ is defined by
$\{x,y\}\mapsto f(x)=f(y)$.
\enddefinition

Here is another easy consequence of Theorem 1.

\proclaim{Corollary} Any odd degree map of a $\Z/2$-homology $n$-sphere to
a stably parallelizable $n$-manifold, $n>2$, is realizable in $\R^{2n}$.
\endproclaim

\demo{Proof} It suffices to consider the case where the given map $f\:N\to M$
is self-transverse.
Let $C$ be a $\tau$-invariant component of $\Sigma_f$.
Consider the diagram
$$\CD C@.\i @.\Sigma_f @>p>> N\\
@VV\Z/2V@. @VV\Z/2V @VVfV\\
C/\tau@.\quad\i\quad @.\Sigma_f/\tau @>f^{(2)}>> M\endCD\tag{$*$}$$
where $p=p_1|_{\Sigma_f}$.
We have $\deg(p|_C)\,\deg(f)=2\,\deg(f^{(2)}|_{C/\tau})$.
Since $\deg(f)$ is odd, $\deg(p|_C)$ must be even, and the assertion follows
from Theorem 1. \qed
\enddemo

\proclaim{Theorem 2} Let $f\:N\to M$ be a self-transverse map between compact
connected oriented $n$-manifolds of non-negative degree, and let
$p_1\:N\x N\to N$ denote the projection.
Then for each component $C$ of $\Sigma_f$, either $0\le\deg(p_1|_C)<\deg(f)$
or $0=\deg(p_1|_C)=\deg(f)$.
\endproclaim

A proof in the case $n=1$ could be based on the standard facts that every
degree $0$ map $S^1\to S^1$ factors through $\R$, whereas every degree $d>0$
map $S^1\to S^1$ factors as a composition of a degree $1$ map and the
$d$-fold cover; the details are left to the reader.
The proof in dimensions $n>1$ is given in \S2.

\proclaim{Corollary} Any degree $0$ map of a $\Z/2$-homology $n$-sphere to
a stably parallelizable $n$-manifold, $n>2$, is realizable in $\R^{2n}$.
\endproclaim

Note that if $\pi_1(M^n)$ is infinite, every map $S^n\to M$ has degree $0$,
factoring through the infinite-sheeted universal covering.

\proclaim{Theorem 3} If $f\:N\to M$ is a self-transverse map between compact
connected oriented $n$-manifolds, where $\pi_1(M)=1$, every component of
$\Sigma_f$ projects to a factor of $N\x N$ with degree $0$.
\endproclaim

The proof, which is given in \S2 for $n>1$, is somewhat in the spirit of
M. Brown's paper \cite{Br}, although no indication of formal connection
has been found.
The proof for $n=1$ is left to the reader.

\proclaim{Corollary} Any map of a $\Z/2$-homology $n$-sphere to a
simply-connected stably parallelizable $n$-manifold, $n>2$, is realizable
in $\R^{2n}$.
\endproclaim

\proclaim{Theorem 4} If $M$ is a closed orientable $3$-manifold, there
exists a map $S^3\to M$, non-realizable in $\R^6$, iff $\pi_1(M)$ is finite
of even order.
\endproclaim

\demo{Proof} If $\pi_1(M)$ is infinite, any map $S^3\to M$ is realizable in
$\R^6$ by Corollary to Theorem 2.
If $\pi_1(M)$ is finite, by Example 3 the universal covering $p\:\hat M\to M$
is such that some $\tau$-invariant component $C$ of $\Sigma_p$ projects with
odd degree to a factor of $\hat M\x\hat M$ iff $|\pi_1(M)|$ is even.

If $|\pi_1(M)|$ is even, choose some $x\in M$ and an arc $J\i\hat M$
containing $p^{-1}(x)$.
By obstruction theory there exists a map $\phi\:S^3\to\hat M$ (which may be
assumed self-transverse) such that each point of $J$ is a regular value of
$\phi$ and has exactly one preimage.
Since $\phi^2$ sends $\Delta_{S^3}$ to $\Delta_{\hat M}$, each component
of $(\phi^2)^{-1}(C)\i\Sigma_{p\phi}$ is compact, hence there is one mapping
onto $C$ with degree $1$.
This component has then to be $\tau$-invariant and project onto a factor of
$S^3\x S^3$ with odd degree.
By Theorem 1, $p\phi$ is non-realizable in $\R^6$.

If $|\pi_1(M)|$ is odd, $p$ is realizable in $\R^6$ by Theorem 1.
By covering theory, any map $f\:S^n\to M$ factors as the composition of
$\hat f\:S^3\to\hat M$ and $p$.
By Corollary to Theorem 3, $\hat f$ is realizable in $\R^6$, hence so is
$f=p\hat f$. \qed
\enddemo

\remark{Remark} Note that the proof of the `if' part works for an arbitrary
$\Z/2$-homology $3$-sphere in place of $S^3$.
Thus, in particular, there exists a self-map of the Poincar\'e homology
sphere, non-realizable in $\R^6$.
\endremark
\medskip

P. M. Akhmetiev noticed recently that an argument due to M. Yamamoto
\cite{Ya} can be used to establish that every map $S^2\to S^2$ realizes in
$\R^4$.
For completeness, we include a proof.

\proclaim{Theorem {\rm (Yamamoto--Akhmetiev)}} Any map $f\:S^2\to M$, where
$M$ is a closed orientable $2$-manifold, is realizable in $\R^4$.
\endproclaim

\demo{Proof} Let us first consider the case $M=S^2$.
Without loss of generality, $f$ is self-transverse, and the poles of $M$ are
regular values of $f$.
Let $p\:M\to I$ be the linear map $S^1*S^0\to B^0*S^0$, where $S^0$ stands
for the poles.
Then $(pf)^{-1}(t)$ is either a finite set (for $t=0,1$), a $1$-manifold
(if $t$ is a regular value of $pf$), or else a $1$-manifold plus one isolated
point or a $1$-dimensional $\Z/4$-manifold with one non-manifold point.
We define a smooth map $g\:S^2\to\R^2$ by sending $(pf)^{-1}(t)$ for each
regular $t$ to a collection of clockwise oriented embedded circles bounding
disjoint disks, and extending by continuity.
Then $(pf)^{-1}(t)$ is embedded by $g$ for each $t\in I$, so
$(pf)\x g\:S^2\to I\x\R^2$ is an embedding.
Hence $f\x g\:S^2\to M\x\R^2\i\R^4$ is an embedding, which is smooth and can
be made arbitrarily close to $f$ by shrinking the neighborhood $M\x\R^2$ of
$M$ in $\R^4$.

Finally, if $M\ne S^2$, $f$ factors as $f\:S^2@>\hat f>>\R^2@>\pi>>M$.
Since $\R^2$ embeds in $S^2$, $\hat f$ is realizable in $\R^4$ by the above.
On the other hand, $\pi$ is realizable in $\R^4$ via the embedding
$\R\to\Gamma\i\R\x M\i\R^4$, where $\Gamma$ denotes the graph of $\pi$.
Thus $f=\pi\hat f$ is also realizable in $\R^4$. \qed
\enddemo

Corollary to Theorem 3 combined with the latter result yield, via
\cite{I; Lemma 2}, a solution to Daverman's problem
\cite{D; Problem E16 (708)}, a context for which can be found in \cite{Mc},
\cite{DH} and \cite{KW}.

\proclaim{Corollary} The limit of any inverse sequence $\{S^n;\,p_i\}$
embeds in $\R^{2n}$ for $n>1$.
\endproclaim

Equivalently, if a compactum $X$ admits for each $\eps>0$ a (continuous)
map $f\:X\to S^n$ such that each point-inverse $f^{-1}(pt)$ is of diameter
$<\eps$, then $X$ embeds in $\R^{2n}$, unless $n=1$.
Solenoids do not embed in $\R^2$ (see \cite{Mc}).

An immediate analog of the double cover $S^1\to S^1$ is the join
$j\:S^3\to S^3$ of its two copies.
Alternatively, $j$ can be described as the composition of the $2$-cover
$S^3\to\RP^3$ and the ramified $2$-fold cover $\RP^3\to S^3$, branched over
the Hopf link.
It is not hard to verify by hand (compare Example 6 in \S0) that some
self-transverse approximation $j'\:S^3\to S^3$ of $j$ has no compact
components in $\Sigma_{j'}$.
However, the `obvious' self-transverse approximation of $j$ in $\R^6$,
obtained by joining linearly the simplest realizations of the double covers
$S^1\to S^1$ in $\R^3\x\{0\}$ and in $\{0\}\x\R^3$, has one double point.
Hence

\proclaim{Conjecture} The join of two dyadic solenoids does not embed in
$\R^6$ in such a way that one of the solenoids goes to $\R^3\x\{0\}$ and
the other to $\{0\}\x\R^3$.
\endproclaim

This paper would never appear without stimulating correspondence and
conversations with Petya Akhmetiev.
I am also grateful to A.~ V.~ Chernavsky, A.~ N.~ Dranishnikov, R.~ Sadykov,
E.~ V.~ Shchepin, A.~ Sz\H ucs and Yu.~ B.~ Rudyak for valuable remarks and
to Juan Liu for help in preparation of the figure.

\head 0. Morin maps \endhead

This section is essentially independent of the rest of the paper, except
that the notation of the following paragraph will be used occasionally.

\definition{Definitions}
Let $f\:N^n\to M^{2n-k}$, $1\le k\le n$, be a self-transverse map between
orientable manifolds; then, as mentioned in the introduction, $\Sigma_f$ is
an orientable submanifold of $N\x N\but\Delta_N$.
Let $\bar\Sigma_f$ stand for the topological closure of $\Sigma_f$ in
$N\x N$.
On the other hand, consider the set $\Delta_f\i N$ of critical points of $f$,
i.e\. those points where $df\:TN\to TM$ is not of the maximal rank.
Identifying $N$ with the diagonal $\Delta_N$ and $TN$ with its normal bundle,
we see that $\bar\Sigma_f\i\Sigma_f\cup\Delta_f$.
It may be verified directly that $\bar\Sigma_f$ is locally euclidean at
the points of $\Delta_f$ of the simplest stable types: simple pinch points
(for $k<n$), fold, cusp, and swallowtail singularities (for $k=n$).
\enddefinition

\proclaim{Proposition 1} If $N$ is compact and $f\:N^n\to M^{2n-k}$,
$1\le k\le n$, is a completely self-transverse Morin map, $\bar\Sigma_f$
is a submanifold of $N\x N$, intersecting $\Delta_N$ transversely in
$\Delta_f$.
\endproclaim

The map $f$ is called a {\it Morin} map, if $\dim(\ker df_x)\le 1$
at every point $x\in N$.
It is not hard to see that the set of Morin maps is open and under
the restriction $2k<n+4$ also dense in the $C^1$ topology (see
\cite{Hi; \S3, proof of Theorem 2.6}).
On the other hand, $f$ is a Morin map if it lifts vertically to an immersion
$N\imm M\x\R$; thus if $k=n$ and $f^*(TM)$ is stably isomorphic to $TN$,
the set of Morin maps is dense in the $C^0$ topology (see \cite{G},
\cite{RS}).

We call the map $f$ {\it completely self-transverse,} if it is
self-transverse and lifts vertically to an immersion
$\bar f\:N\imm M\x\R^k$ such that the monomorphism of the spherical tangent
bundles $S\bar f\:SN\to S(M\x\R^k)$, regarded as a map between the total
spaces, is transverse to the kernel $T_{\R^k}(M\x\R^k)$ of the projection
$T(M\x\R^k)\to TM$.
Note that if $f$ is Morin, the preimage $(S\bar f)^{-1}T_{\R^k}(M\x\R^k)$ is
a pair of sections of the spherical tangent bundle $SN\to N$ over $\Delta_f$,
so that the latter is a manifold whenever $f$ is completely self-transverse.
It is not hard to see that the set of completely self-transverse maps is
open and dense in the $C^1$ topology (cf\. \cite{Mc1}).

\demo{Proof}
Let $\hat N$ denote the blowup of $N\x N$ along $\Delta_N$, which is
a compactification of the deleted product $\tilde N$ by the projective
tangent bundle $\Cal PN$ (i.e\. the projectivization of $SN$).
(In more detail, $\hat N$ is homeomorphic to the manifold obtained from
$N\x N$ by removing the interior of a $\tau$-invariant normal disk bundle
of $\Delta_N$ and identifying points of the remaining sphere bundle via
$\tau$.)
Set $X=\Delta_M\x\tilde{\R^k}\i\widetilde{M\x\R^k}$ and consider
its closure $\bar X:=X\cup\Cal P_{\R^k}(M\x\R^k)$ in the blowup
$\widehat{M\x\R^k}$.
Each $x\in N$ is contained in an $U\i N$ such that $F:=\bar f|_U$ is an
embedding.
Since $f$ is completely self-transverse, the map of the blowups
$\hat F\:\hat U\to\widehat{M\x\R^k}$ is transverse to $\bar X$.
On the other hand, since $f$ is Morin, the canonical projection
$\hat U\to U\x U$ sends $\hat F^{-1}(\bar X)$ homeomorphically onto
$\bar\Sigma_f\cap U\x U$. \qed
\enddemo

By the above remarks, completely self-transverse Morin maps are $C^0$-dense
(and in the case $n<4$ even $C^1$-dense) among all smooth maps between
given stably parallelizable $n$-manifolds.

\remark{Remark {\rm (J. R. Klein and Yu\. B. Rudyak)}} Any (integral)
homology sphere $\Sigma$ is stably parallelizable.
Indeed, the stable tangent bundle of $\Sigma$ is induced
by a map $\phi\:\Sigma\to BO$, whose image is in $BSO$ since $\Sigma$
is orientable.
Since $SO$ is connected, $BSO$ is $1$-connected, and the obstructions
to null-homotopy of $\phi$ lie in the groups $H^i(\Sigma;\pi_i(BSO))$.
The rest of the proof repeats that for homotopy spheres \cite{KM}.
\endremark

\proclaim{Proposition 2} Let $f$ be a completely self-transverse Morin map of
a $\Z/2$-homology $n$-sphere to a stably parallelizable $n$-manifold, $n>2$,
of nonzero even degree.
If $\bar\Sigma_f$ is orientable, $f$ is not realizable in $\R^{2n}$.
\endproclaim

\demo{Proof}
By the hypothesis, any regular value has an even number of preimages under
$f$.
Let $x$ be one, then $x$ has an odd number of preimages under
$p_1|_{\Sigma_f}$.
By Theorem 1, it suffices to show that an odd number of them occur in some
$\tau$-invariant compact component of $\Sigma_f$.
If $C$ is a non-invariant component of $\Sigma_f$, diagram $(*)$ in
introduction implies
$\deg(p_1|_C)\deg(f)=\deg(f^{(2)}|_{C/\tau})=\deg(p_1|_{\tau C})\deg(f)$.
Since $\deg(f)$ is nonzero, $\deg(p_1|_C)=\deg(p_1|_{\tau C})$.
So, if $C$ is compact, $x$ has an even number of preimages under
$p_1|_{C\cup\tau C}$.
If $C$ is a non-compact component of $\Sigma_f$, we enlarge diagram $(*)$
so as to include the limit points of $\Sigma_f$ in the diagonal.
The action of $\Z/2$ is no longer free on the closure $\bar C$ of
$C\cup\tau C$, and the quotient $\bar C/\tau$ has nonempty boundary.
Therefore the composition $(p_1|_{\bar C})^*f^*$ factors through the trivial
group $H^n(\bar C/\tau)$, and since $\deg(f)\ne 0$, the degree of
$p_1|_{\bar C}$ must be zero.
So again $x$ has an even number of preimages under $p_1|_{C\cup\tau C}$,
which completes the proof. \qed
\enddemo

It is shown in \S3 that for every self-transverse $f\:N^n\to M^n$, where $N$
is a $\Z/2$-homology sphere, $M$ is stably parallelizable and $n\ne 1,3,7$,
every $\tau$-invariant component of $\Sigma_f$ projects with even degree to
a factor of $N\x N$.
Thus the proof of Proposition 2 has

\proclaim{Corollary} $\bar\Sigma_f$ is non-orientable for every completely
self-transverse Morin map $f$ of a $\Z/2$-homology $n$-sphere to a stably
parallelizable $n$-manifold, $n\ne 1,3,7$, of nonzero even degree.
\endproclaim

The following lemma implies that $\bar\Sigma_f$ is orientable iff $\Sigma_f$
has no non-compact $\tau$-invariant components.

\proclaim{Lemma 0 {\rm (Akhmetiev)}} If $f\:X\to Y$ is a completely
self-transverse Morin map between orientable $n$-manifolds, the
quotient $\bar\Sigma_f/\tau$ is orientable, meanwhile orientation of
$\bar\Sigma_f$ is reversed along every loop in $\bar\Sigma_f$ crossing
$\Delta_X$ in a single point.
\endproclaim

Notice that $\bar\Sigma_f$ is a two-fold cover of $\bar\Sigma_f/\tau$
branched along the boundary.

\demo{Proof} Let us define a local orientation of $\bar\Sigma_f$ at every
regular point $(p,q)$ of the composition $\chi$ of the inclusion
$\bar\Sigma_f\i X\x X$, the projection $X\x X\to X$ and our $f\:X\to Y$
(note that $\chi$ is the same map no matter onto which of the factors
we project).
We set the orientation at $(p,q)$ to agree with that induced by $\chi$ from $Y$
if and only if the orientations of $X$ and $Y$ either agree (via $f$) at $p$
and agree at $q$, or disagree at $p$ and disagree at $q$.
It is easy to check that these local orientations agree on all of $\Sigma_f$,
but reverse each time we pass through $\Delta_X\cap\bar\Sigma_f$. \qed
\enddemo

\example{Example 5}
Figure 1 depicts a degree two map $f\:S^2\to S^2$, obtained from a
`pseudo-eversion' of $S^1$ in $\R^2$, i.e\. an immersion
$\phi\:S^1\x I\imm\R^2\x I$ with $\phi^{-1}(\R^2\x\{i\})=S^1\x\{i\}$
(`regular pseudo-homotopy') embedding $S^1$ with different orientations
in the top and bottom slices.

\bigskip
\centerline{\epsffile{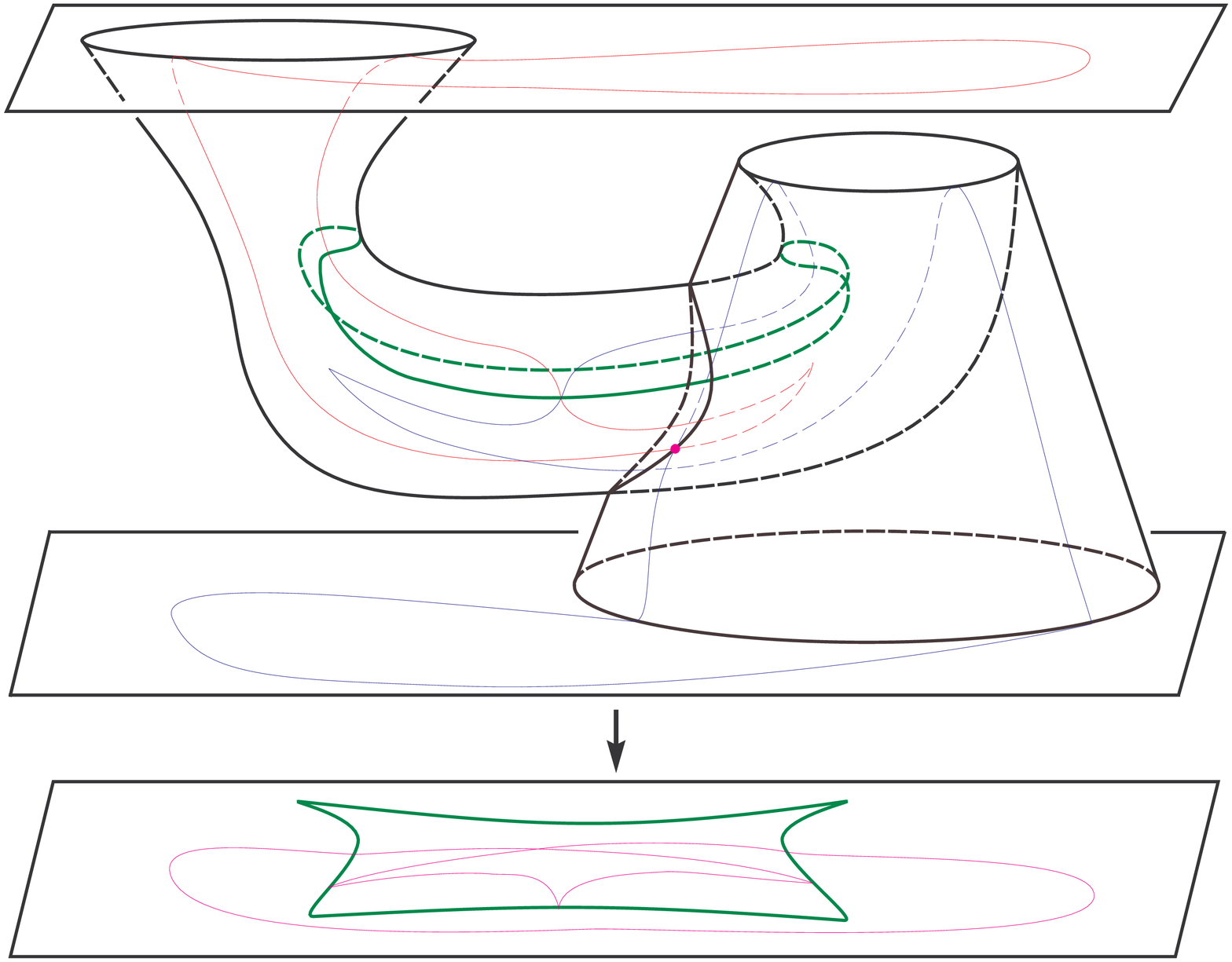}}
\medskip
\centerline{\bf Figure 1}
\bigskip

To see directly that $\bar\Sigma_f$ is non-orientable, we notice that there
is a path $\R\to\Sigma_f$ (shown in Fig\. 1), approaching the same point of
the codimension one submanifold $\bar\Sigma_f\cap\Delta_{S^2}$ from two
different sides in $\bar\Sigma_f$.
By Lemma 0, such a loop is necessarily orientation-reversing.
\endexample

\example{Example 6} The suspension $f\:S^2\to S^2$ of the $2$-cover
$S^1\to S^1$, where $S^1$ is identified with the equator of $S^2$ and
the suspension points with the poles, can be $C^1$-approximated by a map $g$
with $g(\Delta_g)$ consisting of two hyperbolic triangles -- one in
the northern hemisphere and the other in the southern (their edges are the
folds and their vertices are the cusps of $g$; the interior of each triangle
is covered by $4$ sheets).
It might be easier to imagine $g$ by first considering its
$C^0$-approximation $h$ where each northern cusp of $g$ is cancelled with
a southern one by means of a thin tunnel, framed by two folds and going
along a meridian.
This $h$ can be described as the identity on the complement of three disks,
each of which is mapped homeomorphically onto the closure of its complement.

Our goal is to verify by hand that $\Sigma_g$ (as well as $\Sigma_h$) is
non-orientable.
For convenience of notation we assume that $g$ is longitude preserving, that
the polar circles are contained in the interiors of the triangles, and that
each vertex, but no entire edge, is between the tropics.
The restriction of $f$ to the equator $E$ is a $2$-cover, so the intersection
of $\bar\Sigma_g$ with $E\x E$ is a circle (consisting of pairs of antipodal
points), mapped onto $E$ as a $2$-cover by the composition $\chi$ of the
projection $S^2\x S^2\to S^2$ and $f$.
The north tropic $T$ intersects the interior of the northern triangle in
three arcs.
Under $g$, each of these arcs is covered by three sheets in a {\tt Z}-shaped
fashion (resulting in one fold over each end of the arc) and one additional
sheet.
Consequently $\bar\Sigma_g\cap (T\x T)$ consists of $4$ connected
components, and $\chi$ maps one of them onto $T$ with degree two and the others
onto the arcs, so that each arc is covered $6$ times by the first component,
leading to $2$ folds and $2$ regular points at each end of the arc, and
$6$ times by its corresponding component, leading to $3$ folds at each end
of the arc.
Thus there are $5$ folds in total mapping to the eastern end of a fixed arc,
and as we travel further to the north, they will tend to cancel the $5$ folds
mapping to the western end of the arc that is next from the fixed one to
the east.
Since the same must happen along each edge of the northern triangle, our $4$
components have to undergo $15$ surgeries before they can reach the north
polar circle.
But $g$ is the trivial $4$-cover over the north polar circle $P$, so
$\bar\Sigma_g\cap (P\x P)$ is simply $12$ circles, each mapping
homeomorphically onto $N$ under $\chi$.
However, for the reasons of parity, out of fifteen surgeries transforming
four circles into twelve, at least one must be non-orientable. \qed
\endexample

\head 1. Sphere eversions  \endhead

This section is essentially independent of the other ones, though it will
be used in \S2 to give an alternative proof of Theorem 3 in the case of
degree $2$ maps.
Also the proof of Observation 1 below rests on Proposition 1 from \S0, but
this dependence will be eliminated in Proposition 4.

We recall that an {\it eversion} of $S^{n-1}$ is a regular homotopy
$H\:S^{n-1}\x I\imm\R^n\x I$ (where $H^{-1}(\R^n\x\{t\})=S^{n-1}\x\{t\}$ for
each $t\in I$) between a smooth embedding $i\:S^{n-1}\emb\R^n$ and $i$
precomposed with an orientation reversing homeomorphism $\rho$ of $S^{n-1}$.
Specifically, $\rho$ is set to be the antipodal involution if $n$ is odd
(which must be the case if any eversion of $S^{n-1}$ is to exist).
To {\it cap off} an eversion $\phi\:S^{n-1}\x I\to\R^n$ of $S^{n-1}$, we
extend $\phi$ to a degree two map $f_\phi\:S^n\to S^n$, where
$S^n=D_0^n\cup S^n\x I\cup D_1^n$ in the domain and $S^n=\R^n\cup\{\infty\}$
in the range, by identifying $D_0^n$ and $D_1^n$ with the closure $B^n$ of
the component of $S^n\but i(S^{n-1})$ containing the point at infinity.

An eversion $\phi\:S^{n-1}\x I\to\R^n$ is said to be {\it time-symmetric}
if it is equivariant with respect to the involution $(x,t)\mapsto(-x,1-t)$ on
$S^{n-1}\x I$ and the identity map on $\R^n$, and immerses a neighborhood of
$S^{n-1}\x\{\frac12\}$.
Thus $\phi$ factors through the $2$-cover $S^{n-1}\x I\to M$, where $M$ is
the mapping cylinder of the $2$-cover $S^{n-1}\to\RP^{n-1}$, and some map
$M\to\R^n$.
The latter extends to a map $g_\phi\:\RP^n\to S^n$ such that $f_\phi$ is
$g_\phi$ precomposed with the $2$-cover $S^n\to\RP^n$.
Note that although $g_\phi$ may be assumed self-transverse, $f_\phi$ is
never self-transverse, whereas the lift $\bar g_\phi\:\RP^n\imm\R^{n+1}$
of $g_\phi$, defined by the first half of the regular homotopy, has
the same tangent hyperplane $\R^n\x\{\frac12\}$ at each point of
$\bar g_\phi(\RP^{n-1})$.

\proclaim{Observation 1} Let $\phi\:S^{n-1}\x I\to\R^n$ be a time-symmetric
eversion such that $g:=g_\phi$ is a completely self-transverse map.
There exists an arbitrarily close self-transverse $C^1$-approximation
$h\:S^n\to S^n$ of $f:=f_\phi$ such that every $\tau$-invariant component
of $\Sigma_h$ projects with an even degree to a factor of $S^n\x S^n$ if
and only if $g$ maps some essential circle $l\i\RP^n$ onto an arc.
\endproclaim

\demo{Proof} The $2$-cover $\pi\:S^n\to\RP^n$ restricts to a $2$-cover
$\Delta_f\to\Delta_g$ of the fold surfaces.
We choose an $h$ so that it agrees with $f$ outside a neighborhood
$N\Delta_f$ of $\Delta_f$ and on the $\pi$-preimage of at least one point in
each component of $\Delta_g$.
Now $\Sigma_f$ is a union of the antidiagonal $\nabla_{S^n}$ and a (possibly
disconnected) manifold $M$, which intersect transversely along the copy
$\nabla_{\Delta_f}$ of $\Delta_f$.
Thus $\Sigma_f$ is a $\Z/4$-manifold with non-manifold locus
$\nabla_{\Delta_f}$.
The manifold $\Sigma_h$ may be assumed to coincide, up to a small ambient
isotopy, with $\Sigma_f$ outside a neighborhood of $\nabla_{\Delta_f}$.
By our choice of $h$, any two points of
$\nabla_{S^n\but N\Delta_f}\i\Sigma_f$ remain in the same component $N$ of
$\Sigma_h$.
So $N$ is compact iff all the components of $M$ intersecting $\nabla_{S^n}$
are compact.
By considering the point at infinity, $N$, if compact, projects to a factor
of $S^n\x S^n$ with degree one, and no other component of $\Sigma_h$ can
project with nonzero degree.
Thus we only need to determine whether $\nabla_{S^n}$ and $\Delta_{S^n}$ are
contained in distinct connected components of $\Sigma_f\cup\Delta_{S^n}$.

Suppose that they are not.
Notice that $\Sigma_f\cup\Delta_{S^n}$ is a $4$-fold cover of
$\Sigma_g\cup\Delta_{\RP^n}$, namely a restriction of the $4$-cover
$S^n\x S^n\to\RP^n\x\RP^n$.
A path $\ell\:\R\to\Sigma_f\but\nabla_{S^n}$ connecting a point
$\ell(-\infty)=(x,x)\in\Delta_{S^n}$ to a point
$\ell(+\infty)=(y,-y)\in\nabla_{S^n}$, projects to a path
$\ell'\:\R\to\Sigma_g$ connecting the points $(\pi(x),\pi(x))$ and
$(\pi(y),\pi(y))$ of $\Delta_{\RP^n}$.
Let $p_i$ denote the projection onto the $i$-th factor of $S^n\x S^n$, then
$p_1\ell(\R)\cup\{x\}\cup p_2\ell(\R)$ is a path from $y$ to $-y$, so the
union of the projections of $\ell'$ to the factors of $\RP^n\x\RP^n$,
compactified by the endpoints, is a loop not null-homotopic in $\RP^n$.
By construction, the image of this loop under $g$ is an arc.
This completes the proof of the `only if' part; the `if' part follows by
the same arguments. \qed
\enddemo

\example{Example 7}
Using Observation 1, it is particularly easy to see that Shapiro's
time-symmetric eversion \cite{Fra} (without restricting the choice of
the immersion $\RP^2\imm\R^3$ and of the distinguished $\RP^1\i\RP^2$)
yields a map, realizable in $\R^6$.
Indeed, every fiber of the tubular neighborhood of the immersed $\RP^2$,
except for the fibers over the points of the tubular neighborhood of the
distinguished $\RP^1\i\RP^2$, is the image of a nontrivial loop in $\RP^3$.
\endexample

\remark{Remark}
Observation 1 already implies that every time-symmetric eversion $\phi$ of
$S^{n-1}$, $n=3,7$, is regularly homotopic, in the class of time-symmetric
eversions, to a $\psi$ such that $f_{\psi}$ is realizable in $\R^{2n}$.
Indeed, if $h_t\:S^n\to\R^n$ is a homotopy between $i$ and a self-transverse
$C^1$-approximation $j$ of the double cover of some immersion
$\RP^{n-1}\imm\R^n$, define a homotopy $\phi(h_t)\:S^{n-1}\x I\to\R^n$
between $i$ and $i\rho$ by combining $h_t$ and $h_{1-t}\rho$, linked by
the straight line regular homotopy between $j$ and $j\rho$.
It is not hard to construct a non-regular homotopy $h_t$ such that
$g_{\phi(h_t)}$ is self-transverse and sends some essential loop $l$ to
the interior of a generic arc $J$.
By the parametric $C^0$-dense $h$-principle \cite{G}, $h_t$ is arbitrarily
close to a regular homotopy $r_t$ between $i$ and $j$ of a given fiberwise
regular homotopy class.
The preimage of $J$ under a small self-transverse homotopy between
$g_{\phi(h_t)}$ and $g_{\phi(r_t)}$ yields a homology between $l$ and
a component $l'$ of $g_{\phi(r_t)}^{-1}(J)$.
Thus $g_{\phi(r_t)}$ sends the essential loop $l'$ onto an arc.
\endremark
\medskip

\proclaim{Proposition 3} Let $\theta\in H^1(N;\Z/2)$ be a nonzero class, and
$g\:N\to M$ a self-transverse map between compact connected $n$-manifolds
without boundary, $n>2$, such that $g_*\pi_1(N)\i 2\pi_1(M)$.
Any $g$-transverse arc $J\i M$ extends to an embedded circle $\ell\i M$ such
that $g^{-1}(\ell)$ contains an embedded circle $\zeta$, disjoint from
a given arc component $Z$ of $g^{-1}(J)$ and such that $\theta([\zeta])\ne0$.
\endproclaim

This proposition together with Observation 1 is already enough to prove that
for each time-symmetric eversion $\phi$ of $S^2$, the degree two map
$f_\phi\:S^3\to S^3$, obtained by capping off $\phi$ at infinity, is
realizable in $\R^6$ (and similarly for time-symmetric eversions of $S^6$).
Indeed, set $g=g_\phi\:\RP^3\to S^3$, and let $J$ be a small arc through
the point at infinity, consisting entirely of regular values.
Then Proposition 3 with $Z=g^{-1}(J)$ yields a circle $\zeta\i\RP^3$, which
does not lift to $S^3$ and is mapped by $g$ onto an arc.

The idea of the proof is due to P. M. Akhmetiev.

\demo{Proof} Let $\eta$ be any embedded circle in $N$, containing the arc $Z$
but otherwise generic with respect to $g$, and such that
$\theta([\eta])\ne 0$.
Thus $\eta$ is transverse to the fold surface $\Delta_f$, and since $n>2$, we
may assume that $\ell:=g(\eta)$ is an embedded circle in $M$.
By the hypothesis, $[\ell]$ is a $\bmod 2$ boundary.
Then so is $[g^{-1}(\ell)]=g_![\ell]$, where the transfer $g_!$ is
the composition of the Poincar\'e duality in $M$, the cohomological
$g^*$ and the Poincar\'e duality in $N$.
Now $g^{-1}(\ell)$ is the union of $\eta$ and a closed $1$-manifold $\eta'$,
where $\eta'\cap\eta$ is precisely the set of non-manifold points of
$g^{-1}(\ell)$.
Note, however, that $g_![\ell]=[\eta]+[\eta']$ only modulo $2$ (cf\. Remark
to Lemma 2 in \S2).
Since $\theta([\eta'])=\theta([\eta])\ne 0$, there exists a connected
component $\zeta$ of $\eta'$ such that $\theta([\zeta])\ne 0$.
But $\zeta$ may meet $\eta$ only at the non-manifold points of
$g^{-1}(\ell)$, which do not occur in the arc $Z$ of $\eta$, since $J$ is
transverse to $g$. \qed
\enddemo

\proclaim{Lemma 1} Let $f\:I\sqcup I'\to I$ be any map such that
$f(0\sqcup 0')=0$ and $f(1\sqcup 1')=1$.
Then $(0,0')$ and $(1,1')$ belong to the same component of $\Sigma_f$.
\endproclaim

\demo{Proof} We may assume that $0$ and $1$ are regular values of $f$
(otherwise extend $f$ by the trivial $2$-cover over a collar).
Then both $(0,0')$ and $(1,1')$ belong to the component of
$\Sigma_f\cap I\x I'$ that projects to $I$ with degree $1$ rel $\partial$.
\qed
\enddemo

The same argument as for time-symmetric eversions above, with Observation 1
replaced by the following proposition, proves that every degree two map
$S^3\to S^3$, obtained by capping off at infinity any (not necessarily
time-symmetric) eversion of $S^2$ in $\R^3$, is realizable in $\R^6$ (and
similarly for eversions of $S^6$).

\proclaim{Proposition 4} Let $f\:X\to Y$ be any map and $M$ a
$\tau$-invariant connected component of $\Sigma_f$.
If there exists a loop $l\:S^1\to M/\tau$ that does not lift to $M$ and
is such that $f^{(2)}l$ factors through $\R$, then $M$ is non-compact.
\endproclaim

\demo{Proof} Let $g$ denote the map $f^{(2)}\:\Sigma_f/\tau\to Y$; we are
going to mimic the proof of Observation 1.
By Lemma 1, there exists a path
$\ell\:\R\to\Sigma_g$, connecting two points of $\Delta_{M/\tau}$, say,
$\ell(-\infty)=(\{a,b\},\{a,b\})$ and $\ell(+\infty)=(\{c,d\},\{c,d\})$, and
such that the union of the projections of $\ell(\R)$ to the factors of
$M/\tau\x M/\tau$, compactified by the endpoints, does not lift to $M$.
This lifts to a path $\bar\ell\:\R\to\Sigma_{g\pi}$, where $\pi\:M\to M/\tau$
denotes the double cover, connecting $((a,b),(a,b))$ to either $((c,d),(d,c))$
or $((d,c),(c,d))$ --- say, to the first one.
Let $P_1$, $P_2$ denote the projections onto the factors of $M\x M$; then
$\ell_1:=P_1\bar\ell\:\R\to M$ connects $(a,b)$ and $(c,d)$, meanwhile
$\ell_2:=P_2\bar\ell$ connects $(a,b)$ and $(d,c)$, moreover the compositions
of these two paths with $g\pi=fp_1=fp_2$ coincide, where $p_i$ denote
the projections onto the factors of $X\x X$.
Hence $\ell_{11}:=p_1\ell_1\:\R\to X$ connects $a$ to $c$, meanwhile
$\ell_{22}:=p_2\ell_2\:\R\to X$ connects $b$ to $c$, and of course the
compositions of these paths with $f$ coincide.
Thus $\ell_{11}\x\ell_{22}\:\R\to\Sigma_f\i X\x X$ connects $(a,b)\in M$
to $(c,c)\in\Delta_X$. \qed
\enddemo

\proclaim{Corollary} If $M$ is a $1$-manifold and $f\:M\to S^1$
a self-transverse map, any component of $\Sigma_f/\tau$, mapped by $f^{(2)}$
through an arc, lifts to $\Sigma_f$.
\endproclaim

The assumption that $f$ is self-transverse cannot be omitted, as shown by
the composition $f$ of the $2$-cover $S^1\to S^1$ and any degree $0$ map
$S^1\to S^1$.
(Incidentally, such a composition does not lift vertically to an embedding
$S^1\emb S^1\x\R$.)

\demo{Proof} Suppose that $C$ is a component of $\Sigma_f$ that does not
lift to $\Sigma_f$ and is mapped by $f^{(2)}$ through an arc.
By Proposition 4, $C$ is non-compact.
Then $C$ is contractible and therefore does lift to $\Sigma_f$. \qed
\enddemo

\head 2. Unfolding self-transverse maps \endhead

This section contains proofs of Theorems 2 and 3.
It is essentially independent of the rest of the paper, except for using
Lemma 1 and the proof of Lemma 0 from preceding sections.

\proclaim{Lemma 2} Let $f\:N\to M$ be a self-transverse map between
connected $n$-mani\-folds, $n>2$.
If $x,y\in N$ are $f$-regular points, any $f$-transverse arc between $f(x)$
and $f(y)$ extends to an $f$-transverse circle $\ell\i M$ such that $x$ and
$y$ lie in the same component $C$ of $f^{-1}(\ell)$.
\endproclaim

This lemma is already enough to prove Theorem 2 in the case where some
$f$-regular value $z\in M$ has exactly $2$ preimages $x$, $y$.
Indeed, the restriction of $f$ to the circle $C$ is a degree $0$ map
$\phi\:S^1\to S^1$ with $\phi^{-1}(z)=\{x,y\}$, and by Lemma~ 1, $(x,y)$
belongs to a non-compact component of $\Sigma_\phi$.

\demo{Proof}
Let $J\i M$ be the given arc between $f(x)$ and $f(y)$, and let $p\:I\to N$
be a path connecting $x$ to $y$ and generic with respect to $f$.
Thus $p$ is transverse to the fold surface $\Delta_f$ and, since $n>2$, we
may assume that $\ell':=J\cup fp(I)$ is a smoothly embedded circle in $M$.
However, $\ell'$ will not in general be $f$-transverse, for $\ell'$ has to be
tangent with $f(\Delta_f)$ at the $f$-image of any point of
$p(I)\cap\Delta_f$.
So $\Gamma:=f^{-1}(\ell')$ is a $1$-dimensional $\Z/4$-manifold (i.e\. a
graph where each vertex has valency $4$), and since it contains $p(I)$,
the points $x$ and $y$ belong to the same component $\Gamma_0$ of $\Gamma$.
If we slightly perturb $\ell'$ to remove its tangencies with $f(\Delta_f)$,
$\Gamma$ will resolve into a genuine $1$-manifold, with $\Gamma_0$ possibly
becoming disconnected.
The perturbation can be done in two ways at each tangency point, which may
either vanish or turn into two generic intersection points.
If, with some choice of the perturbation, $\Gamma_0$ resolves to $k>1$
circles, there must be two of them adjacent to the same vertex of $\Gamma_0$.
Changing the type of perturbation at this vertex connects the two circles
into one.
Thus the perturbation of $\ell'$ can be chosen so that $\Gamma_0$ resolves
to a single circle. \qed
\enddemo

\remark{Remark}
The combinatorial part of the above argument is nothing but a proof of
existence of an oriented Eulerian circuit in the oriented graph $\Gamma_0$,
where each vertex has incoming valency $2$ and outgoing valency $2$.
(An orientation on the edges of $\Gamma$ is induced by an orientation of
$\ell'$ via $f_!:H_1(\ell';\Z/4)\to H_1(\Gamma;\Z/4)$, which is the
composition of the Poincar\'e duality, $f^*$ and the Poincar\'e duality.)
\endremark

\definition{Terminology} Let $f\:N\to M$ be a self-transverse map between
compact oriented $n$-manifolds, and let $J\i M$ be an $f$-transverse arc.
Any orientation of $J$ yields a co-orientation of $J$ in $M$, which via $f$
induces a co-orientation of $f^{-1}(J)$ in $N$, hence an orientation of
$f^{-1}(J)$.
A component $C$ of $f^{-1}(J)$ is called a {\it positive (negative)} arc if
$f|_C\:(C,\partial C)\to (J,\partial J)$ has degree $+1$ (resp\. $-1$).
(The sign depends on the orientations of $M$ and $N$, but not on that of $J$.)
Else $C$ could be a circle or an arc with both endpoints mapping onto
the same endpoint of $J$, with $f(C)\ne J$.
\enddefinition

\proclaim{Lemma 3} Let $f\:S^1\to S^1$ be a self-transverse map, $J$ an arc
in the range whose endpoints $x$, $y$ are regular values, and suppose that
$J_+$ is a positive and $J_-$ a negative arc in $f^{-1}(J)$.
Then there exists an embedded path $p\:(I,\partial I)\emb (S^1,f^{-1}(x))$,
either starting with $J_+$ and ending with some negative arc $J'_-$ or
starting with $J_-$ and ending with some positive arc $J'_+$, such that
$fp\:(I,\partial I)\to (S^1,\{x\})$ factors through $(\R,\{0\})$.
\endproclaim

\demo{Proof} By symmetry we may assume $\deg(f)\ge 0$.
Then there are at least as many positive arcs in $f^{-1}(J)$ as negative.
Starting from $J_-$ and moving along $S^1$, we will at some point (before
completing the full circle) reach a positive arc $J_+'$ such that on our way
from $J_-$ to $J_+'$ we will have passed as many positive arcs in $f^{-1}(J)$
as negative (excluding both $J_-$ and $J_+'$).
If $p$ denotes our path, $fp\:(I,\partial I)\to (S^1,\{x\})$ has degree $0$
and therefore factors through $(\R,\{0\})$. \qed
\enddemo

\remark{Remark} It is easy to see that the constructed path factors through
$(\R_+,\{0\})$, where $\R_+$ denotes the non-negative reals, provided that
we stop at the very first occurrence of a suitable arc $J_+'$, resp\. $J_-'$.
\endremark

\proclaim{Observation 2} Let $f\:N\to M$ be a self-transverse map between
compact connected oriented $n$-manifolds, $n>2$.
If $\deg(f)\ge 0$, there exists an $f$-transverse arc $\Cal J\i M$ such that
$f^{-1}(\Cal J)$ contains no negative arcs.
\endproclaim

Equivalently, the algebraic sum of the $f$-preimages of an $f$-regular value
$z\in\Cal J$ (with their $f$-signs) in any component of $f^{-1}(\Cal J)$ must
be non-negative.

\demo{Proof} Let $J\i M$ be an arbitrary $f$-transverse arc with an endpoint
$x$.
Let $m_J$ denote the number of negative arcs in $f^{-1}(J)$.
It suffices to find, assuming that $m_J>0$, an arc $J'$ such that
$m_{J'}<m_J$.
Let $J_+$ and $J_-$ be a positive and a negative arc in $f^{-1}(J)$.
By Lemma 2, $J$ extends to an $f$-transverse circle $\ell\i M$ such that
$J_+$ and $J_-$ are contained in the same connected component $C$ of
$f^{-1}(\ell)$.
By Lemma 3, there is an embedded arc $p(I,\partial I)\i (C,f^{-1}(x))$,
containing either $J_+$ and some negative arc $J'_-$ or $J_-$ and some
positive arc $J'_+$, and such that $fp\:(I,\partial I)\to (\ell,\{x\})$
factors as $(I,\partial I)@>\phi>>(\R,\{0\})@>\psi>>(\ell,\{x\})$.
By compactness, $\phi(I)\i [-r,r]$ for some $r\in\N$; without loss of
generality, $\psi([0,1])=J$.
Let $\psi'\:[-r,r]\to M$ be an $f$-transverse embedding, closely
$C^1$-approximating $\psi|_{[-r,r]}$ and agreeing with it on $[0,1]$.
Then the arc $J':=\psi'([-r,r])$ contains $J$, so every negative arc of
$f^{-1}(J')$ contains at least one negative arc of $f^{-1}(J)$, thus
$m_{J'}\le m_J$.
Moreover, $J_+$ and $J_-'$ (resp\. $J_+'$ and $J_-$) are now contained in
the same connected component of $f^{-1}(J')$, which is either not a negative
arc, or contains at least one more negative arc of $f^{-1}(J)$.
Thus $m_{J'}<m_J$. \qed
\enddemo

\remark{Remark} The final part of the proof can be replaced by an algebraic
argument, using the following reformulation of the definition of $m_J$:
$$2m_J+\deg(f)=\sum_{C\in\pi_0(f^{-1}(J))}\left|\sum_{q\in C\cap f^{-1}(z)}
\sign_f(q)\right|.$$
(Here $\pi_0$ stands for the set of connected components.)
\endremark

\remark{Remark} Given an $f$-transverse arc $J\i M$ with an endpoint $x$,
the arc $\Cal J$ can be chosen to contain $J\but\{x\}$ as an open subset
(using Remark to Lemma 3).
\endremark

\remark{Remark} Alternatively, given a regular value $z\in M$, the arc
$\Cal J$ can be chosen to contain $z$ and so that each component of
$f^{-1}(\Cal J)$, meeting $f^{-1}(z)$, is either a positive arc or a circle.
Indeed, the original arc given by Observation 2, which we here denote by $J$,
may be assumed to contain $z$ in its interior.
Let $x$ be an endpoint of $J$, and let $J_x$ denote the subarc of $J$ from
$z$ to $x$.
Each arc in $f^{-1}(J)$ with both endpoints in $f^{-1}(x)$ and meeting
$f^{-1}(z)$ corresponds to one negative arc in $f^{-1}(J_x)$.
Apply the previous remark to $J_x$ to get an arc $J'$ with endpoints $z$
and $x'$ such that $f^{-1}(J')$ contains no negative arcs.
Then $f^{-1}(J'\cup J)$ contains no arcs with both endpoints in $f^{-1}(x')$
and meeting $f^{-1}(z)$.
Similarly, we may find an arc $J''$, meeting $J$ in $\Cl{J\but J'}$, so that
$\Cal J:=J'\cup J\cup J''$ is as required.
\endremark

\proclaim{Lemma 4} Let $f\:N\to M$ be a self-transverse map between compact
connected oriented $n$-manifolds.
Suppose that $\deg(f)\ge 0$ and there exists an $f$-transverse arc $J\i M$
such that $f^{-1}(J)$ contains precisely $\deg(f)$ positive arcs $J_i$.
If $x$ is an endpoint of $J$ and $x_i$ the endpoint of $J_i$ in $f^{-1}(x)$,
every component $L$ of $\Sigma_f$ projects to a factor of $N\x N$ with degree
equal to the number of pairs $(x_1,x_i)$ contained in $L$.
\endproclaim

\demo{Proof} Let $x_+$ and $x_-$ be the endpoints of some component $C$ of
$f^{-1}(J)$ such that $f(x_+)=f(x_-)=x$.
If $\deg(f)>0$, the points $(x_1,x_+)$ and $(x_1,x_-)$ are contained in
the same component of $\Sigma_f$ by Lemma 1 and have opposite signs with
respect to the projection to the first factor of $N\x N$ (cf\. proof of
Lemma 0).
Thus the preimage of $x_1$ in $L$ may contain some of the points
$(x_1,x_i)$ and some number of couples $(x_1,x_+)$, $(x_1,x_-)$, where
$x_+$ and $x_-$ are as above, and each such couple contributes zero to
the algebraic sum of the signs of the preimages of $x_1$ in $L$.

If $\deg(f)=0$, consider the component $\hat C$ of $f^{-1}(J)$ such that
any component $C$ as above satisfies $f(C)\i f(\hat C)$.
Let $\hat x_+$, $\hat x_-$ be the endpoints of $\hat C$ and $x_+$, $x_-$
those of some $C$ as above.
By Lemma 1 the component of $\Sigma_f$ containing $(\hat x_+,\hat x_-)$ is
non-compact and $(\hat x_+,x_+)$, $(\hat x_+,x_-)$ are contained in
the same component of $\Sigma_f$ due to our choice of $\hat C$.
Thus $L$ projects with degree $0$ to the first factor. \qed
\enddemo

\demo{Proof of Theorem 2} If $n>2$, by Observation 2 there exists an
$f$-transverse arc $J\i M$ such that $f^{-1}(J)$ contains precisely $\deg(f)$
positive arcs, and the assertion follows from Lemma 4.
For $n=2$ the proof is analogous, replacing all arcs and circles by immersed
arcs and circles. \qed
\enddemo

\remark{Remark} The material of the preceding section can now be used
to prove Theorem 3 in the case of degree $2$ maps.
(This will not be used in the sequel.)
Assume that $n>2$; when $n=2$, embedded arcs will have to be replaced by
immersed.
Let $f\:N\to M$ be a self-transverse degree $2$ map, and $C$ a
$\tau$-invariant component of $\Sigma_f$, projecting with nonzero degree onto
a factor of $N\x N$.
By Observation 2 there exists an $f$-transverse arc $J\i M$ such that
$f^{-1}(J)$ contains precisely $2$ positive arcs $J_1$, $J_2$.
By Lemma 4, $C$ contains the component $\hat Z$ of $\Sigma_f\cap J_1\x J_2$,
projecting with degree $1$ onto $J_1$.
Let $\eta$ be the line bundle associated with the $2$-cover $C\to C/\tau$,
and $w_1(\eta)\in H^1(C/\tau)$ its first Stiefel--Whitney class.
Since $\pi_1(M)=1$, by the proof of Proposition 3 with $\theta=w_1(\eta)$ and
$Z=(\hat Z\cup\tau\hat Z)/\tau$ there exist a circle $\ell\i M$, containing
$J$, and a loop $\zeta\:S^1\to (C/\tau\cap(f^{(2)})^{-1}(\ell))\but Z$, which
does not lift to $C$.
Indeed, although $f^{(2)}$ is not necessarily self-transverse, it is locally
self-transverse (i.e\. becomes self-transverse when restricted to
a sufficiently small neighborhood of every point of the domain), and
the proof of Proposition 3 goes through, except that $\eta'$ will now be
a $\Z/4$-manifold, and consequently the resulting circle $\zeta$ immersed
rather than embedded.
(An alternative approach would be to note that $f^{(2)}$ is homotopic by
a small homotopy through locally self-transverse maps to a self-transverse
map.)
Since $J_1$ and $J_2$ are the only positive arcs in $f^{-1}(J)$,
$Z$ is the only component of $(f^{(2)})^{-1}(J)$ whose image covers $J$,
and therefore $f^{(2)}\zeta\:S^1\to M$ factors through an arc.
By Proposition 4, $C$ is non-compact, which is a contradiction.
\endremark

\proclaim{Lemma 5} Let $f\:(M^2,\partial M)\to (D^2,\partial D)$ be
a self-transverse map of a compact connected oriented $2$-manifold onto
the disk, $\deg(f)\ge 0$, and suppose that there exists an $f$-transverse arc
$J\i\partial D$ such that $f^{-1}(J)$ contains precisely $\deg(f)$
positive arcs $J_i$.
Let $x$ be an endpoint of $J$ and $x_i$ the endpoints of $J_i$ in $f^{-1}(x)$.
For any pair $x_i\ne x_j$, contained in the same component of $\partial M$,
the point $(x_i,x_j)$ is contained in a non-compact component of $\Sigma_f$.
\endproclaim

See Example 6 in \S0 for a model situation.
The assertion of Lemma 5 does not hold if the disk $D^2$ is replaced with
the punctured torus $H^2$, by considering the cover of $H$ corresponding to
any homomorphism of $\pi_1(H)=\Z*\Z$ onto a non-abelian finite group.

\demo{Proof} Let $\phi_t\:S^1\to D^2$ be a self-transverse isotopy from
the inclusion $\phi_0$ onto $\partial D^2$ to an embedding $\phi_1$ into
a small neighborhood of $J$, such that $\phi_t|_J=\id$ for each $t$.
Let $f_t\:X_t\to S^1$ be the restriction of $f$ to $f^{-1}(\phi(S^1))$,
where each $X_t$ is a $1$-manifold except for a finite number of values
of $t$, for which it is a $\Z/4$-manifold with one non-manifold point $x_t$.
Then $f_0=f|_{\partial M}$, each $X_t$ contains $f^{-1}(J)$, and each
component of $X_1$ maps with degree $+1$ or $0$ onto $\phi_1(S_1)$ according
as it contains a positive arc of $f^{-1}(J)$ or not.
Rename the positive arcs $J_i$, $J_j$ as $J_+$, $J_-$ and their endpoints
$x_i$, $x_j$ as $x_+$, $x_-$.
Then $x_+$ and $x_-$ are in the same component of $\partial M=X_0$,
but in different components of $X_1$.
Let $s\in (0,1)$ be the smallest critical level of the bordism $f_t$ such
that $x_+$ and $x_-$ are in the same component $\ell_{-\eps}$ of $X_{s-\eps}$
but in distinct components $\ell^+_{+\eps}$, $\ell^-_{+\eps}$ of
$X_{s+\eps}$.
The canonical projection $r_{+\eps}\:X_{s+\eps}\to X_s$ (with precisely one
non-degenerate point-inverse) takes the latter
onto oriented submanifolds $r_{+\eps}(\ell^+_{+\eps})$ and
$r_{+\eps}(\ell^-_{+\eps})$ of the $\Z/4$-manifold $X_s$, whose union is
$\ell_0:=r_{-\eps}(\ell_{-\eps})$ and whose intersection is the non-manifold
point $\{x_s\}$.

Let us reverse the orientation of $r_{+\eps}(\ell^-_{+\eps})\i X_s$ and
consider the corresponding oriented resolution $X_{s+i\eps}$ and its
canonical projection $r_{i\eps}$.
Then $\ell_{i\eps}=r_{i\eps}^{-1}(\ell_0)$ is connected.
With respect to the map
$f_{s+i\eps}:=f_sr_{i\eps}\:X_{s+i\eps}\to\phi_s(S^1)$, the arc $J_+$ is
again positive, but $J_-$ is now negative in $f_{s+i\eps}^{-1}(J)$.
Applying Lemma 3 to the restriction of $f_{s+i\eps}$ to $\ell_{i\eps}$,
we find a path $p\:I\to\ell_{i\eps}@>r_{i\eps}>>\ell_0$ that either
starts with $J_+$ and ends with some negative arc $J'_-$ or starts with $J_-$
and ends with some positive arc $J'_+$, moreover
$f_sp\:(I,\partial I)\to (S^1,\{x\})$ lifts to the universal cover
$(\R,\{0\})$.
On the other hand, we can certainly find immersed paths
$p_\pm\:I\to\ell^\pm_{+\eps}@>r_{+\eps}>>\ell_0$, starting at $J_\pm$ and
such that the image of the lift of $f_sp_\pm\:(I,\{0\})\to(S^1,\{x\})$ to
$(\R,\{0\})$ contains that of $f_sp$.
By Lemma 1, applied to $p\sqcup p_-$ (resp\. $p\sqcup p_+$), the point
$(x_+,x_-)$ lies in the same component of $\bar\Sigma_{f_s}$ with
$(x'_-,x_-)$ (resp\. $(x_+,x'_+)$), where $x'_\pm$ denotes the endpoint
of $J'_\pm$ in $f^{-1}(x)$.

If $x'_-=x_-$ (resp\. $x'_+=x_+$), we are done; otherwise $x'_-$ and $x_-$
(resp\. $x'_+$ and $x_+$) are in the same component $\ell^-_{+\eps}$
(resp\. $\ell^+_{+\eps}$) of $X_{s+\eps}$, but in distinct components of
$X_1$, and the proof proceeds by induction on the number of critical levels.
\qed
\enddemo

\demo{Proof of Theorem 3} First, let us assume $n>2$; for some orientations
of $M$ and $N$ we have $\deg(f)\ge 0$.
By Observation 2 there exists an $f$-transverse arc $J\i M$ such that
$f^{-1}(J)$ contains no negative arcs.
Suppose that $C$ is a component of $\Sigma_f$ such that $\deg(p_1|_C)\ne 0$.
By Lemma 4, there exist positive arcs $J_1$, $J_2$ in $f^{-1}(J)$ with
endpoints $x_1$, $x_2$ such that $(x_1,x_2)\in C$.
By Lemma 2, $J$ extends to an $f$-transverse circle $\ell$ such that
$J_1$ and $J_2$ are contained in the same component of $f^{-1}(\ell)$.
Since $\pi_1(M)=1$, by the relative $C^0$-dense $h$-principle \cite{G},
\cite{RS} $\ell$ bounds an immersed disk in $M$, which can of course be
assumed $f$-transverse.
By Lemma 5, the component $C$ is non-compact, which is a contradiction.

The proof for $n=2$ is the same, except that all embedded arcs and circles
must be replaced by immersed ones, regularly null-homotopic in the case of
$\ell$. \qed
\enddemo

\head 3. The controlled stable Hopf invariant \endhead

This section contains proofs of Theorem 1 and of (explicit and implicit)
results of Akhmetiev \cite{A1}, \cite{A2} mentioned in the introduction.
It is entirely independent of the preceding sections.
\smallskip

Consider again a self-transverse map $f\:X^n\to Y^{n+k}$ between orientable
manifolds, and recall the notation $\tl X=X\x X\but\Delta$.
Let $\Z\T$ denote the local coefficient system on $\tl X/\tau$, associated
with the $2$-cover $\tl X\to\tl X/\tau$; more precisely, it could be
described either as the $0$-dimensional Leray sheaf of this cover, or as
determined by the $\Z[\Z/2]$-module $\Z[\Z/2]/(\tau+1)$ with respect to
the action of $\Z/2=\left<\tau\mid2\tau\right>$ on $\tl X$.
We recall that (co)homology of the quotient $P/\pi$ of a polyhedron $P$ by
a free PL action of a finite group $\pi$ with coefficients in the local
system $\Cal O_M$ determined by a $\Z\pi$-module $M$ is just the homology of
the chain complex $C_*\otimes_{\Z\pi}M$ (resp\. $\Hom_{\Z\pi}(C_*;M)$),
where $C_*$ denotes the complex of integral chains of a $\pi$-invariant
triangulation of $P$.
Similarly, locally finite homology $H_*^\lf(P/\pi;\Cal O_M)$ (with possibly
infinite cycles) and cohomology with compact support $H^*_c(P/\pi;\Cal O_M)$
can be defined as the homology groups of the chain complexes
$\Hom_{\Z\pi}(C^*_c;M)$
and $C^*_c\otimes_{\Z\pi}M$, respectively, where $C^*_c$ denotes the
subcomplex of cochains with finite support in the integral cochain complex
$C^*=\Hom(C_*,\Z)$.
Since $\tau$ preserves the orientation of $\tl X$ iff $n$ is even, and the
co-orientation of $\Sigma_f$ in $\tl X$ (which is induced from the
co-orientation of $\Delta_Y$ in $Y\x Y$) iff $n+k$ is even, the orienting
sheaf of $\Sigma_f/\tau$ is $\Z\T^{\otimes k}$ (where the tensor product can
be understood in the sense of either sheaves or modules), which is $\Z\T$ if
$k$ is odd, and $\Z$ if $k$ is even.
In other words, there is the element of infinite order
$[\Sigma_f/\tau]\in H_d^\lf(\Sigma_f/\tau;\Z\T^{\otimes k})$, where $d=n-k$.

\proclaim{Lemma 6} {\rm (compare \cite{Mi; Lemma 1, Theorem 1} and
\cite{Mc2; proof of Theorem 4.6})} Let
$f\:N^n\to M^{n+k}$ and $\bar f\:N\to M\x\R$ be self-transverse maps between
oriented manifolds such that $f=\pi\bar f$, where $\pi$ denotes
the projection, and let $\eta$ denote the line bundle associated with
the $2$-cover $\Sigma_f\to\Sigma_f/\tau$; set $d=n-k$.

(a) $[\Sigma_{\bar f}/\tau]\in H_{d-1}^\lf(\Sigma_f/\tau;\Z/2)$
is Poincar\'e dual to $w_1(\eta)\in H^1(\Sigma_f;\Z/2)$;

(b) $[\Sigma_{\bar f}/\tau]\in
H_{d-1}^\lf(\Sigma_{\!f}/\tau;\Z\T^{\otimes(k+1)})$
is dual to the Euler class $e(\eta)\in H^1(\Sigma_{\!f}/\tau;\Z\T)$.
\endproclaim

Of course, $e(\eta)$ turns into $w_1(\eta)$ once we reduce the coefficients
$\bmod 2$, which does not lead to a loss of information
since $e(\eta)$ anyway is an element of order two.
However, the assertion (b) is strictly stronger than (a), for it implies that
$[\Sigma_{\bar f}]$ has order two in the (twisted) integral homology.

\demo{Algebraic proof} Define $\phi\:\Sigma_f\to(\Sigma_f/\tau)\x\R$ by
$(x,y)\mapsto(\{x,y\},\bar f(x)-\bar f(y))$.
Notice that $\phi^{-1}((\Sigma_f/\tau)\x\{0\})=\Sigma_{\bar f}$ and
$p\phi=\pi$, where $\pi\:\Sigma_f\to\Sigma_f/\tau$ is the $2$-cover and $p$
projects $(\Sigma_f/\tau)\x\R$ onto the first factor.
Hence $\phi$ defines a bundle morphism $\Lambda$ between $\eta$ and
the trivial line bundle $\eps$ over $\Sigma_f/\tau$, moreover
$\Sigma_{\bar f}/\tau$ is the locus of points where $\Lambda$ fails to be
injective.
So the adjoint $\Lambda^*$ yields a section of $\eta$ with zero set
$\Sigma_{\bar f}/\tau$. \qed
\enddemo

\demo{Geometric proof} The element of $[\Sigma_f/\tau,\RP^d]$ classifying
the $2$-cover $\Sigma_f\to\Sigma_f/\tau$ can be represented by the map
$\phi$, assigning to every unordered pair $\{x,y\}$ the line through
the points $\hat f(x)$ and $\hat f(y)$, where
$\hat f\:N^n\to M^{n+k}\x\R^{n-k+1}$ is an embedding (which exists by general
position) projecting to $f$ along the second factor.
If moreover $\hat f$ projects to $\bar f$ along the last $n-k=d$ coordinates,
the homotopy class of the restriction of $\phi$ to $\Sigma_{\bar f}/\tau$
(whose image is contained in $\RP^{d-1}$) classifies the $2$-cover
$\Sigma_{\bar f}\to\Sigma_{\bar f}/\tau$.
Hence $\Sigma_{\bar f}/\tau$ can be regarded as a transversal preimage of
$\RP^{d-1}$ under $\phi$.
Since $\RP^{d-1}$ is the zero set of a section of the universal line bundle
$\gamma$ over $\RP^d$, $\Sigma_{\bar f}/\tau$ is the zero set of the induced
section of $\eta=\phi^*\gamma$. \qed
\enddemo

\remark{Remarks} (i) As a motivation for the algebraic argument, which is
essentially contained in \cite{AS}, we notice that $\eps$ is the pullback of
$\nu_{M\x\R}(M)$ under $f^{(2)}$, meanwhile
$(f^{(2)}\x\id_\R)\o\phi\:\Sigma_f\to M\x\R$ can be identified with
the composition of $\bar f^2|_{\Sigma_f}$ and the projection
$\Delta_M\x\R\x\R\to\Delta_M\x\nabla_{\R}$, where $\nabla_{\R}$ denotes
the antidiagonal $\{(t,-t)\mid t\in\R\}$.
(ii) The geometric argument, which follows the philosophy of \cite{KS}, can
be made somewhat more invariant, if, instead of choosing a section of
$\gamma$, we notice that, by transversality, the normal bundle of
$\Sigma_{\bar f}/\tau$ in $\Sigma_f/\tau$ is induced from that of $\RP^{d-1}$
in $\RP^d$, and apply naturality of Thom class.
\endremark

\proclaim{Lemma 7} Let $f\:N^n\to M^{n+k}$ and $g\:N\to M\x\R^i$ be
self-transverse maps between oriented manifolds such that $f$ is
$\eps$-close to $if$, where $i$ denotes the inclusion of $M=M\x\{0\}$.
Let $\bar O_\eps\i\tl N$ denote the $\eps$-neighborhood of
$\Sigma_f\cup\Delta_N$ in $N\x N$, and set
$O_\eps=(\bar O_\eps\but\Delta_N)/\tau$.
Let $w_1$ and $e$ denote the first Stiefel--Whitney class and the Euler class
of the line bundle associated with the $2$-cover $\tl N\to\tl N/\tau$;
set $d=n-k$.

(a) $[\Sigma_g/\tau]=[\Sigma_f/\tau]\!\Cap\! w_1^i$, \
$\Cap\:H_d^\lf(O_\eps;\Z/2)\otimes H^i(\tl N/\tau;\Z/2)\to
H_{d-i}^\lf(O_\eps;\Z/2)$;

(b) $[\Sigma_g/\tau]=[\Sigma_f/\tau]\!\Cap\! e^i\!$, \
$\Cap\:H_d^\lf(O_\eps;\Z\T^{\otimes k})\otimes
H^i(\tl N\!/\?\tau;\Z\T^{\otimes i})\to
H_{d-i}^\lf(O_\eps;\Z\T^{\otimes(k+i)})$.
\endproclaim

\demo{Proof} It suffices to consider the case $i=1$.
Let $\bar f\:N\to M\x\R$ be a self-transverse map, $\eps$-close to $g$ and
such that $f=p\bar f$, where $p$ denotes the projection.
Then the assertion with $g$ replaced by $\bar f$ follows from Lemma 6.
Let $H\:N\x I\to M\x\R\x I$ be a self-transverse $\eps$-homotopy between
$\bar f$ and $g$.
Then $\Sigma_H/\tau\i O_\eps\x\Delta_I$ yields an embedded bordism between
$\Sigma_{\bar f}/\tau$ and $\Sigma_g/\tau$.
In particular, the homology classes in question coincide. \qed
\enddemo

\proclaim{Theorem 5} Let $f\:N^n\to M^{n+k}$ be a self-transverse map between
orientable manifolds, where $n>2$ and $N$ is compact, and let
$i\:M=M\x\{0\}\emb M\x\R^{n-k}$ denote the inclusion.
Then $if$ is $\eps$-approximable by embeddings for every $\eps>0$ iff
$w_1^{n-k}(\eta)=0$, where $\eta$ denotes the line bundle associated with
the $2$-cover $\Sigma_f\to\Sigma_f/\tau$.
\endproclaim

Equivalently, $w_1^{n-k}(\eta)[C]=0\in\Z/2$ for every connected component
$C$ of $\Sigma_f/\tau$.

\remark{Remark}
The stable Hopf invariant also featured in the context of equivalence of
embeddings $X^n\emb\R^{2n}$ whose images are contained in
$\R^{2n-k}\i\R^{2n}$ \cite{Gi}.
\endremark

\remark{Remark} If $(X,\tau)$ is a space with involution, the maximal $d$
such that $w_1^d\ne 0$, where $w_1$ denotes the first Stiefel--Whitney class
of the line bundle associated with the $2$-cover $X\to X/\tau$, is known as
the Yang index of $(X,\tau)$.
It gives an upper bound for the maximal $d$ such that there exists
an equivariant map $(S^d,t)\to (X,\tau)$, where $t$ denotes the antipodal
involution, and a lower bound for the minimal $d$ such that there exists
an equivariant map $(X,\tau)\to (S^d,t)$, see \cite{CF}.
In this connection, note that an embedding $g\:N\emb M\x\R^{n-k}$ such that
$p_1g=f$ yields an equivariant map $(\Sigma_f,\tau)\to (S^{n-k},t)$,
defined by $(x,y)\mapsto\frac{p_2g(x)-p_2g(y)}{||p_2g(x)-p_2g(y)||}$; here
$p_i$ denote the projections onto the factors of $M\x\R^{n-k}$.
\endremark

\remark{Remark} It is not hard to prove, using Proposition 1, the proof
of Theorem 5 below and some arguments from \cite{A3}, that a $C^1$-generic
map $f\:N^n\to M^{2n-k}$, which is $C^0$-approximable by embeddings in
$M\x\R^{k-l}$, lifts vertically to an embedding $\bar f\:N\emb M\x\R^{k-l}$,
provided that $2(l+1)<k\le n$ and $l\le 2(n-k)$ (compare
\cite{A3; Conjecture 1.9}).
The details will appear elsewhere.
\endremark

\demo{Proof}
Suppose there exists an embedding $S^n\emb M\x\R^{n-k}$, $\eps$-close to $if$
for a sufficiently small $\eps>0$ (defined below).
By Lemma 7(a), $[\Sigma_f/\tau]\Cap w_1^{n-k}=0$ in $H_0^\lf(O_\eps;\Z/2)$.
Any homeomorphism between $\tl N/\tau$ and $(N\x N\but DN)/\tau$, where $DN$
is a $\tau$-equivariant tubular neighborhood of $\Delta_N$, yields a
compactification $\tl N/\tau\cup PN$ of $\tl N/\tau$, where the corona
$PN\cong\partial(DN)/\tau$ can be identified with the total space of
the projective tangent bundle of $N$.
This leads to the excision isomorphism
$H_i^\lf(O_\eps;\Z/2)\simeq H_i(O_\eps\cup PN,PN;\Z/2)$.
If $\eps>0$ is small enough, there is a deformation retraction of
$O_\eps\cup PN$ onto $\Sigma_f/\tau\cup PN$ fixing $PN$.
So the latter group is isomorphic with
$H_i(\Sigma_f/\tau\cup PN,PN;\Z/2)\simeq H_i^\lf(\Sigma_f/\tau;\Z/2)$.
It is easy to trace that the equation $[\Sigma_f/\tau]\Cap w_1^{n-k}=0$
survives through these isomorphisms to the very last group.
Now $[\Sigma_f/\tau]$ becomes the fundamental class in
$H_{n-k}^\lf(\Sigma_f/\tau;\Z/2)$, so we get
$w_1^{n-k}(\eta)=0\in H^{n-k}(\Sigma_f/\tau;\Z/2)$.

Conversely, suppose the latter holds, then
$e^{n-k}(\eta)=0\in H^{n-k}(\Sigma_f/\tau;\Z\T^{\otimes k})$ by the remark
following Lemma 6.
Consider a self-transverse map $\bar f\:N^n\to M\x\R^{n-k}$ such that
$f=p\bar f$, where $p$ denotes the projection.
By the iterated application of Lemma 6b (as in the proof of Lemma 7),
$[\Sigma_{\bar f}/\tau]=0$ in the homology of $\Sigma_f/\tau$, hence in that
of $O_\eps$ for every $\eps>0$ (in the notation of Lemma 7).
But this (twisted) integral homology class is easily identified as
the Poincar\'e dual of the van Kampen--Skopenkov obstruction to
$\eps$-approximability by an embedding, which is complete for $n>2$ \cite{Sk}
(see also \cite{AM}).
Since $\eps>0$ could be taken arbitrarily small, $f$ is discretely
realizable. \qed
\enddemo

Since any embedding of a stably parallelizable $n$-manifold in $\R^{2n}$ has
trivial normal bundle (see introduction), Theorem 1 is a special case of

\proclaim{Corollary} Let $f\:N\to M$ be a self-transverse map of
a $\Z/2$-homology $n$-sphere to an orientable $n$-manifold, $n>2$, and let
$i\:M\i M\x\R^n$ denote the inclusion.
Then $if$ is $\eps$-approximable by embeddings for every $\eps>0$ iff
every $\tau$-invariant connected component of $\Sigma_f$ projects with an even
degree to a factor of $N\x N$.
\endproclaim

\demo{Proof} Let $w_1$ denote the first Stiefel--Whitney class of the line
bundle associated with the $2$-cover $\tl N\to\tl N/\tau$.
By the K\"unneth formula and the exact sequence of pair, combined with
the Thom isomorphism, $\tl N$ has the $\bmod 2$ cohomology of $S^n$.
By the Smith sequence, the $\bmod 2$ cohomology ring of $\tl N/\tau$ is
the same as that of $\RP^n$.
So if $C/\tau$ is a compact component of $\Sigma_f/\tau$, $w_1^n[C/\tau]$ is
nonzero iff $[C/\tau]\ne 0$ in $H_n(\tl N/\tau;\Z/2)=\Z/2$.
By the homological Smith sequence, the latter is equivalent to
$[C]\ne 0$ in $H_n(\tl N;\Z/2)=\Z/2$.
Since $\tau$ acts trivially on this group, $C$ must be $\tau$-invariant
if it is to be not null-homologous $\bmod 2$.
Finally, since $i_*\:H_n(\tl N;\Z/2)\to H_n(N\x N;\Z/2)$ is the monomorphism
onto the diagonal, $[C]\ne 0$ iff $C$ projects with odd degree to a factor of
$N\x N$. \qed
\enddemo

Geometrically, when $N=S^n$, for every compact component $C/\tau$ of
$\Sigma_f/\tau$ the $\bmod 2$ residue $w_1^n\Cap [C/\tau]$ counts the parity
of the number of intersections between $C/\tau$ and $(T_p\cup\tau T_p)/\tau$,
where $T_p=\{p\}\x(S^n\but\{p\})$ for a fixed $f$-regular point $p\in S^n$.
This observation is useful for a parametric version of controlled embedding
theory (cf\. \cite{AM}), so we include the details.

\proclaim{Lemma 8} There exists an equivariant homotopy equivalence
$\phi\:\Tl S^n\to S^n$, trans\-versal to $S^0$ and such that
$\phi^{-1}(S^0)=T_p\cup\tau T_p$, where $T_p=\{p\}\x(S^n\but\{p\})$.
\endproclaim

\demo{Proof}
We will exhibit an equivariant deformation retraction of
$(\Tl S^n,T_p\cup\tau T_p)$ onto $(\nabla_{S^n},\{(p,-p),(-p,p)\})$.
For $(x,y)\in\Tl S^n\but\nabla_{S^n}$ let $S^1_{xy}$ denote the large circle
of $S^n$ containing the geodesic $J_{xy}$ between $x$ and $y$, and $p_{xy}$
the point of this circle closest to $p$, if it is unique.
If it is not unique, $x$ and $y$ will run by $S^1_{xy}$ away from $J_{xy}$
with equal constant speed, until they become antipodal.
Otherwise, $x$ and $y$ will run along $S^1_{xy}$ away from $J_{xy}$ until they
become antipodal, with constant speeds $v_x$ and $v_y$, defined as follows.
The ratio $\frac{v_a}{v_b}=d+(1-d)\frac{l_a}{l_b}$, where $a=x$, $b=y$ if
$l_x\le l_y$ and $a=y$, $b=x$ otherwise, $d=\frac2\pi\dist(p,p_{xy})$,
and $l_x$ is the length of the arc in $S^1_{xy}\but\{y\}$ connecting $x$ to
$p_{xy}$ (if $y=p_{xy}$, the choice of the arc does not matter), and
similarly $l_y$ is the length of the arc in $S^1_{xy}\but\{x\}$ connecting
$y$ to $p_{xy}$. \qed
\enddemo

For completeness, we deduce Akmetiev's result, mentioned in the introduction,
from Theorem 5.

\proclaim{Theorem {\rm (Akhmetiev \cite{A2})}} If $n\ne 1,2,3,7$, any map
$f\:N\to M$ between stably parallelizable $n$-manifolds, where $N$ is
compact, is realizable in $\R^{2n}$.
\endproclaim

\demo{Proof} Since $M$ and $N$ are stably parallelizable, by the $C^0$-dense
$h$-principle for immersions \cite{G} (see also \cite{RS; 4.4}),
the composition of $f$ and the inclusion $M\i M\x\R$ can be approximated by
an immersion $g$.
We will prove that the composition of $g$ and the inclusion $M\x\R\i\R^{2n}$,
given by any section of the trivial normal bundle of the given embedding
$M\emb\R^{2n}$, is approximable by embeddings.
By Theorem 5, it suffices to show that $w_1^{n-1}(\eta)=0$, where $\eta$
denotes the line bundle associated with the cover $\Sigma_g\to\Sigma_g/\tau$.
The normal bundle of $g$ is orientable, hence trivial, so the normal bundle
$\nu$ of $g^{(2)}\:\Sigma_g/\tau\imm M\x\R$ is isomorphic to
$\eta\oplus\eps$, where $\eps$ denotes the trivial line bundle
\cite{KS; p.~202} (see also \cite{Mi; proof of Lemma 2}).
Hence $w_1(\eta)=w_1(\nu)=w_1(\Sigma_g/\tau)$.
By the Hirsch lemma \cite{RS}, $M\x\R$ immerses in $\R^{n+1}$, and
$\Sigma_g/\tau$, immersed in $M\x\R\imm\R^{n+1}$, can be compressed into
$\R^n$; so it remains to apply the following lemma. \qed
\enddemo

\proclaim{Lemma 9} $w_1^{n-1}(K^{n-1})$ vanishes for every manifold $K^{n-1}$
immersible in $\R^n$, provided $n\ne 1,3,7$.
\endproclaim

A direct geometric proof of Lemma 9 (in terms of immersions)
is given in \cite{AS} in the case where $n+1$ is not a power of $2$ and in
\cite{A4} in the case $n\equiv 7\pmod 8$, $n\ne 15$.
This yields a geometric proof, with the exception of dimension 15, of
the Adams theorem on the vanishing of the stable Hopf invariant.

\demo{First proof} By \cite{B} (see also \cite{Li}), if $K^{n-1}$ immerses in
$\R^n$, it is either cobordant to $\RP^0$, $\RP^2$ or $\RP^6$, or
null-cobordant.
In the latter case all its Stiefel--Whitney numbers vanish. \qed
\enddemo

\demo{Second proof} It suffices to consider the case where $K$ is connected.
Applying Koschorke's figure $8$ construction \cite{K} (see also \cite{KS}) to
a self-transverse immersion $\phi\:K\imm\R^n$, we get an oriented immersion
$\psi\:L^n\imm\R^n\x\R$ of the total space $L$ of the orientation circle
bundle over $K$ (i.e\. the Whitney join of the orienting and the trivial
$2$-covers over $K$).
A small perturbation of the last coordinate will make $\psi$ into a
self-transverse immersion $\chi$.
Notice that $\Sigma_\chi/\tau$ is the disjoint union of $K$ and an additional
manifold $K'$; moreover, the composition of $\chi^{(2)}$ and the projection
$\R^n\x\R\to\R^n$ restricts to $K$ as $\phi$ and takes $K'$ into
a neighborhood of $\phi^{(2)}(\Sigma_\phi/\tau)$.
On the other hand, let $\bar\phi\:K\imm\R^{2n-1}$ be a generic embedding
such that $p\bar\phi=\phi$, where $p$ denotes the projection along
the first $n-1$ coordinates.
Then the Koschorke construction $\hat\psi\:L\imm\R^{2n-1}\x\R$ on $\bar\phi$
satisfies $p\hat\psi=\psi$; moreover, we can apply the same perturbation of
the last coordinate to get an immersion $\hat\chi\:L\imm\R^{2n-1}\x\R$ such
that $p\hat\chi=\chi$.
The advantage of $\hat\chi$ over $\chi$ is that it has no additional double
points, that is, $\hat\chi^{(2)}\:\Sigma_{\hat\chi}/\tau\to\R^{2n}$ coincides
with the composition of $\bar\phi\:K\imm\R^{2n-1}$ and the inclusion
$\R^{2n-1}\i\R^{2n}$.
Finally, by a small deformation of $\hat\chi$ in the first $n-1$ coordinates
we get a self-transverse immersion $\bar\chi\:L^n\imm\R^{2n}$ such that
$p\bar\chi=\chi$.
Then $\Sigma_{\bar\chi}/\tau\i\Sigma_{\hat\chi}/\tau=K$, meanwhile Lemma 6(a)
implies that $[\Sigma_{\bar\chi}/\tau]$ is Poincar\'e dual to
$w_1^{n-1}(\Sigma_\chi/\tau)$ (using that, by \cite{KS; p.~202},
$w_1(\Sigma_\chi/\tau)=w_1(\zeta)$, where $\zeta$ denotes the line
bundle associated with the $2$-cover $\Sigma_\chi\to\Sigma_\chi/\tau$),
so actually $[\Sigma_{\bar\chi}/\tau]$ is Poincar\'e dual to $w_1^{n-1}(K)$.
But $[\Sigma_{\bar\chi}/\tau]\in H^0(K;\Z/2)=\Z/2$ is the stable Hopf
invariant of the element of the stable stem, represented by the (framed)
immersion $\chi$ \cite{KS; p.~203}, so by the Adams theorem it must be
trivial whenever $n\ne 1,3,7$. \qed
\enddemo

\remark{Remark}
It follows from the above proof of Akhmetiev's theorem that for any
self-transverse $f\:S^n\to M^n$, where $M$ is stably parallelizable,
$n\ne 1,3,7$, every $\tau$-invariant component of $\Sigma_f$ projects with
even degree to a factor of $S^n\x S^n$.
Indeed, $w_1^{n-1}\Cap [\Sigma_{\bar f}/\tau]=0\in
H_0(\Sigma_{\bar f}/\tau;\Z/2)$ by the proof of Akhmetiev's theorem, whereas
by Lemma 6, $[\Sigma_{\bar f}/\tau]=w_1\Cap[\Sigma_f/\tau]\in
H_{n-1}(\Sigma_f/\tau;\Z/2)$.
Hence $w_1^n\Cap [\Sigma_f/\tau]=0\in H_0(\Sigma_f/\tau;\Z/2)$, and
the assertion follows as in the proof of Corollary to Theorem ~5.
\endremark
\medskip

Akhmetiev's Theorem combined with \cite{I; Lemma 2} yield

\proclaim{Corollary {\rm (Akhmetiev)}} If $n\ne 1,2,3,7$, the limit of any
inverse sequence of compact stably parallelizable $n$-manifolds embeds in
$\R^{2n}$.
\endproclaim

\proclaim{Problem} Does there exist an inverse sequence of compact stably
parallelizable $2$-, $3$- or $7$-manifolds whose limit does not embed in
$\R^4$, $\R^6$, resp\. $\R^{14}$?
\endproclaim

Note that the product of $n$ solenoids, and more generally every inverse
limit of tori $(S^1)^n$ embeds in $\R^{2n}$ for $n=2,3,7$ \cite{KW}.
On the other hand, gluing together the trivial (disconnected) and
the non-trivial $2$-covers of $\RP^n\but D^n$ yields the unique irregular
$2$-cover $\RP^n\#\RP^n\to\RP^n\#\RP^n$, which is non-realizable in
$\R^{2n}$, $n=3,7$ \cite{DH; \S10}.

\enddocument
\end